\documentclass[12pt,reqno]{amsart}
\usepackage{amsmath}
 \usepackage{graphics} 
  \usepackage{epsfig} 
\usepackage{graphicx}  \usepackage{epstopdf}
 \usepackage[colorlinks=true]{hyperref}
\hypersetup{urlcolor=blue, citecolor=red}

 \numberwithin{equation}{section}    

\newtheorem{theorem}{Theorem}[section]
\newtheorem{corollary}[theorem]{Corollary}

\newtheorem{lemma}[theorem]{Lemma}

\theoremstyle{definition}
\newtheorem{definition}[theorem]{Definition}

\newcommand{\ep}{\varepsilon}
\newcommand{\eps}[1]{{#1}_{\varepsilon}}
\newcommand{\set}{\Omega}

\newcommand{\bset}{{\partial\Omega}}
 \newcommand{\R}{{\bf R}}
\newcommand{\dn}{{D_{\nu}}}
\newcommand{\dnz}{{D_{\nu_z}}}
 \newcommand{\ov}{\overline}
 \newcommand{\sv}{\mathcal{V}}
 \newcommand{\sm}{\mathcal{M}}
 \newcommand{\scy}{\mathcal{C}}
   \newcommand{\sd}{\mathcal{D}}
   \newcommand{\sw}{\mathcal{W}}  
 
  \newcommand{\cy}{\mathcal{C}}
 \newcommand{\db}{d_\partial}
 \newcommand{\se}{\mathcal{E}}

    \newcommand{\ch}{{\raisebox{.15 em}{$\chi$}}}
     
\title[Green's function for the Lapacean] 
      {Cancellation Properties and Pointwise Bounds for the Green's Functions for the Laplace Operator}

\author[David Hoff]{}


 \email{hoff@iu.edu}

\begin{document}
\maketitle

\centerline{\scshape David Hoff}
\medskip
{\footnotesize
 \centerline{Indiana University}
   \centerline{Department of Mathematics}
   \centerline{Bloomington, IN, USA}
} 
\today
%
%
\bigskip

 \bigskip

We derive a cancellation property satisfied by the derivatives of the \break Green's functions for the Laplace operator corresponding to Dirichlet and Neumann boundary conditions on bounded sets in $\R^n$. The main result is derived in a broader, self-contained exposition which includes construction of and basic pointwise bounds for the Green's functions and their derivatives. We also give an application of the cancellation property to a problem in fluid mechanics.
\medskip

\section{Introduction}

We derive a cancellation property satisfied by derivatives of the Green's functions $G_D$ and $G_N$ for the Laplace operator corresponding respectively to Dirichlet and Neumann boundary conditions on subsets $\set$ of $\R^n$, $n\ge 3$. This cancellation property is expressed in terms of pointwise bounds independent of distances to $\bset$ and generalizes Newton's third law concerning equal and opposite forces. Bounds for arbitrary derivatives $D_x^\alpha D_y^\beta G_{D,N}$ are also derived, partial results being known but not otherwise easily accessed in the literature. The analysis applies certain ``interior-type estimates at the boundary" which although well established in greater generality and corresponding technicality, admit fairly straightforward derivations for second-order scalar elliptic equations; these derivations are therefore included as well. Besides introducing the cancellation property, which is our main result, we thus give  a complete and self-contained development of the important pointwise bounds satisfied by the Green's functions. We also describe the particular problem in fluid mechanics whose resolution requires the cancellation property and which motivated these results.

More specifically we extend to $G_D$ and $G_N$ versions of the following properties of the Newtonian potential $\Gamma$, which for $n\ge 3$ and distinct points $x$ and $y$ in $\R^n$ is given by 

\begin{equation}\label{1.0}\Gamma(x,y)=-c_n^{-1}|x-y|^{2-n}
\end{equation}
where $c_n$ is the area of the unit sphere in $\R^n$. 
One easily checks that for multi-indices $\alpha$ and $\beta$ there is a constant $C$ depending on $|\alpha|$ and $|\beta|$ such that
\begin{equation}\label{1.1}
|D_x^{\alpha}D_y^\beta \Gamma(x,y)|\le C|x-y|^{2-n-(|\alpha|+|\beta|)}
\end{equation}
and that
\begin{equation}\label{1.2}
\nabla_x \Gamma(x,y) + \nabla_y \Gamma (x,y) =0,
\end{equation}
where $\nabla_x$ and $\nabla_y$ are the usual gradients with respect to $x$ and $y$. The cancellation in \eqref{1.2} expresses the familiar fact concerning equal and opposite gravitational and electrostatic forces in $\R^3$. We will show that $G_D$ and $G_N$ satisfy \eqref{1.1} for $|\alpha|$ and $|\beta|$ in a range determined by the regularity of $\bset$ and that while the exact cancellation \eqref{1.2} does not hold in general, a particular such property does:
\begin{equation}\label{1.3}
|\big(\sigma(x)\cdot\nabla_x + \sigma(y)\cdot\nabla_y\big)D_x^\alpha D_y^\beta G_{D,N}(x,y)|\le C|x-y|^{2-n-(|\alpha|+|\beta|)}
\end{equation}
for $x$ and $y$ near $\bset$ and for local vector fields $\sigma$ which are tangential on $\bset$. Observe that the rate of blowup as $y\to x$ in \eqref{1.3} is \emph{one order less} than that given by \eqref{1.1} for either of the gradients on the left. This improvement is crucial for the application of interest in which the term on the left in \eqref{1.3} is required to be locally integrable for $|\alpha|+|\beta|=1$, even though second derivatives of $G_{D,N}$ are not.

It is instructive to look at two examples. First, for the lower half-space in $\R^n$, $G_D$ and $G_N$ are known to be $\Gamma-\Gamma_{rfl}$ and $\Gamma+\Gamma_{rfl}$ respectively, where $\Gamma_{rfl}(x,y)=\Gamma(x,\ov y)$ and $\ov y$ is the reflected point $(y_1,\ldots,y_{n-1},-y_n)$. Direct computation shows that 
$$(\partial/\partial x_j)G_{D,N}+(\partial/\partial y_j)G_{D,N}=0$$
 for $j< n$ but not for $j=n$, consistent with \eqref{1.3} but also showing that \eqref{1.3} need not hold for nontangential vector fields $\sigma$.

A more interesting example is that in which $\set$ is the unit ball in $\R^n$ centered at the origin. We let $G$ denote either $G_D$ or $G_N$ and exploit the symmetry by considering a rotation $Q(\theta)$ through angle $\theta$ about an arbitrary axis through the origin. Anticipating basic facts about Green's functions (which are derived in detail in Section~3 and with which most readers will already be familiar), we compute by changing variables that the Green's function for $Q\set$ is $G(Q(\theta)x,Q(\theta)y)$. But $Q\set=\set$ and Green's functions are unique, so that  $G(Q(\theta)x,Q(\theta)y)=G(x,y)$. Taking the derivative with respect to $\theta$ at $\theta=0$ and observing that $\dot Q(0)x\equiv x^\perp$ is perpendicular to the outer normal vector $x\in\bset$ and similarly for $y$, we thus find that 
$$\big(x^\perp\cdot\nabla_x + y^\perp\cdot\nabla_y\big)G(x,y)=0$$
as in \eqref{1.2}. This is a stronger result than in \eqref{1.3}, and this same argument can be applied more generally to achieve exact cancellation when the set $\set$ exhibits rotational symmetries, as for tori and ellipsoids. 
Exact cancellation as in \eqref{1.2} for  general domains $\set$ is not anticipated however, and is likely inaccessible by our methods. 

Our analysis consists primarily in approximation by Greens functions for half-spaces determined by planes tangent to $\bset$ together with elementary scaling estimates for standard imbedding theorems applied on small sets. The approximation arguments are quite technical because the approximation operator must be understood in sufficient detail to allow for its commutation with the differential operator on the left in \eqref{1.3}. The scaling estimates are of note primarily for their usefulness, the proofs being fairly simple. As we shall see, the pointwise bounds of interest for $G_{D,N}(x,y)$ are relatively straightforward when either $x$ or $y$ is far from $\bset$ or when $x$ and $y$ are far from each other. The main difficulty therefore occurs when $x$ and $y$ are close and both are close to $\bset$, and in this case we require the ``interior-type estimates near the boundary" for weak solutions of general elliptic second-order equations referred to above. These show that for small balls $B'$ and $B$ centered at a point on $\bset$ with $B'$ compactly contained in $B$, a higher-order norm of a solution of such an equation on $B'\cap\set$ can be estimated in terms of its lower-order norm on $B\cap\set$. These estimates are given for both Lebesgue and Sobolev norms (but on cylinders rather than balls in order to more easily accommodate our preference for representing $\bset$ locally in terms of graphs rather than level sets).

The plan of the paper is as follows. First, after collecting various definitions and notations in Section~2, we give a complete, elementary derivation in Section~3 of the basic properties of the Green's functions $G_D$ and $G_N$, the main results given in Theorem~3.1. In Section 4 we state and prove pointwise bounds for $D_x^\alpha D_y^\beta G_{D,N}(x,y)$ similar to those in \eqref{1.1} detailing rates of blowup as $y\to x$, bounds which are independent of distances to $\bset$, together with results concerning H\"older continuity of derivatives. These are given in Theorem 4.1, and in Corollary~4.2 we derive more basic bounds which do depend on distances to $\bset$ and which are required in subsequent sections. The main result of the paper is Theorem~5.1 in which we prove \eqref{1.3} together with two useful variations. In Section 6 we show how the bounds in Theorems~4.1 and 5.1 are applied in resolving the instantaneous regularization of the velocity field for viscous, compressible fluid flow in sets $\set\subset \R^3$. Finally in the Appendix we state and prove in Theorems~7.2 and 7.3 interior-type estimates near the boundary for general elliptic equations \and in Theorem ~7.4 the scaling results referred to above.  
\medskip

There is a considerable literature on Green's functions for the Dirichlet problem for general elliptic equations and systems; see for example \cite{gruter}, \cite{zhen}, \cite{conlon}, \cite{hwang}, \cite{mourgoglu} and the references therein. This literature includes some but not all the results in Theorem 4.1 for the Laplacian with Dirichlet boundary conditions. For the Neumann case the subject is somewhat less developed however: \cite{kenig} and \cite{dahlberg} include results tangentially related to those of the present paper, but concerning pointwise bounds we find only \cite{robert} and \cite{hoff3d}. Theorem 4.1 for the Neumann case therefore extends the literature in a significant way, and our main result, Theorem 5.1, is new. Definitive interior-type estimates at the boundary, including for higher order equations and systems, are derived in the classical papers \cite{agmon1} and \cite{agmon2}.\bigskip

\section{Preliminaries}

We begin with some preliminary facts and notations. First, the ball centered at $x$ in $\R^n$ or $\R^{n-1}$ is denoted $B_R(x)$, or simply $B_R$ if $x=0$; and $\cy_R$ is the cylinder $B_R\times (-R,R)\subset \R^n$. The spatial domain $\set$ will be fixed throughout and is described as follows:
\smallskip

\begin{definition} \emph{$\set$ will be bounded, open and connected in $\R^n$ with a $C^{k+1}$ boundary, where $n\ge 3$ and $k\ge 0$. This means that there are positive constants $R_0$ and $M_\set(k)$ such that for every $z\in \bset$ there is a rigid motion $T^z$ of $\R^n$ (that is, a translation followed by a rotation) and a $C^{k+1}$ mapping $\psi^z$ such that for all $R\in (0,R_0]$, $\psi^z$ maps $B_R\subset\R^{n-1}$ into $(-R/2,R/2)$ and the following hold:
\begin{itemize}
\item $T^z(z)=0$\smallskip
\item $\psi^z(0)=0$ and $\nabla\psi^z(0)=0$\smallskip
\item the derivatives of $\psi^z$ up to order $k+1$ are bounded in absolute value by $M_\set(k)$ in $B_{R_0}$\smallskip
\item $T^z(\set)\cap \mathcal C_{R_0}= \{(\tilde y,y_n) : \tilde y \in \  B_{R_0}\  {\rm and}\ -R_0<y_n<\psi^z(\tilde y)\}$.
\end{itemize}
It will be convenient to define the inverse mapping $S^{z}=(T^z)^{-1}$ and sets $\scy^z_R=S^z(\scy_R)$ and} 
\begin{equation}\label{2.1}
\sw^{z}_R=S^z \big(\{(\tilde y,y_n) : \tilde y \in \  B_{R}\  {\rm and}\ -R<y_n<\psi^{z}(\tilde y)\}\big).
\end{equation}
\end{definition}
\smallskip

\smallskip\noindent \noindent With these assumptions there is a $C^k$ vector field $\nu$ defined in a neighborhood of $\bset$ whose restriction to $\bset$ is the unit outer normal vector field  on $\bset$, a nonnegative Radon measure $dS$ representing surface area on the Borel sets in $\bset$ and for $p\in[1,\infty)$ a bounded linear transformation from $W^{1,p}(\set)$ into $L^p(\bset)$ whose action on a continuous function in $W^{1,p}(\set)$ is its restriction to $\bset$. ($W^{j,p}(\set)$ is the Banach space of elements of $L^p(\set)$ having weak derivatives up to order $j$ which also are in $L^p$, and $W^{2,j}\equiv H^j)$; and the divergence theorem holds for integrands in $W^{1,p}(\set$). Finally, the set of restrictions to $\set$ of $C^\infty$ functions on $\R^n$ is dense in $W^{j,p}(\set)$ for $p\in[1,\infty)$, this fact applying as well to open sets satisfying a ``segment condition." See \cite{hoffbook} sections 3.2 and A.4 and \cite{adams} Theorem 3.22 for complete statements and proofs. 

The outward normal derivative of a sufficiently smooth function $u$ on $\ov\set$ is the directional derivative $\nabla u\cdot\nu$, where $\nabla$ is the usual gradient; this normal derivative will be denoted $\dn u$. Observe that if $u\in W^{2,p}(\set)$ then $\dn u \in L^p(\bset)$ by the remarks in the preceding paragraph.

 We denote the distance from a point $y\in\ov\set$ to $\bset$ by $d_\partial(y)$. The following fact is easily checked: {\rm There is a positive number $\kappa_0<R_0$ depending on $M_\set(0)$ such that if
$d_\partial(y)\le\kappa_0$ then
then there is a unique $y'\in\bset$ such that $d_\partial(y)=|y'-y|.$} 

We denote by $C(\set,k,\dots)$ a generic positive constant depending only on $n,\,k$ and the parameters $R_0 $ and $M_\set(k)$ in {\rm Definition 2.1}, on the diameter of $\set$ and on other parameters listed in the ellipsis, if any. The value of $C(\set,k,\dots)$ may change from step to step but is allowed only to increase as the proofs proceed. Similar conventions will hold for other constants $R_1,\,R_2$ and $R_3$ which will appear in later sections. 

The usual Laplace operator is denoted by $\Delta$ and the average value of an element $u\in L^1(\set)$, that is, its integral over $\set$ divided by the measure of $\set$, will be denoted by $\ov u$. The diagonal in $\R^n\times \R^n$ is the set $\{(x,x) : x\in\R^n\}$ which will be denoted by $\sd$, and for $\ep>0$, $\sd_\ep\equiv\{(x,y)\in \R^n \times \R^n : |x-y|\le \ep\}$.

 A mapping $u$ from a subset $\sw$ of $\R^n$ into $\R$ is H\"older continuous with exponent $\lambda\in (0,1]$ if the seminorm
$$\langle u\rangle_{\sw,\lambda}\equiv \sup\frac{|u(x)-u(y)|}{|x-y|^\lambda},$$
where the sup is over distinct points $x$ and $y$ in $\sw$, is finite. If $\sw$ is open, $C^{j,\lambda}(\sw)$ is the set of functions whose derivatives up to order $j$ exist and are bounded and H\"older-$\lambda$ in $\sw$, so that the norm
$$|u|_{C^{j,\lambda}(\sw)}\equiv\sum_{|\alpha|\le j}\big({\rm sup}_\sw|D^\alpha u|+\langle D^\alpha u\rangle_{\sw,\lambda}\big)$$
is finite.
Note that if  $u\in C^{j,\lambda}  (\sw)$ then $u$ and its derivatives up to order $j$ have continuous extensions to $\partial\sw$. 

\bigskip

\section{Construction of the Green's Functions} First we review some basic facts concerning the Newtonian potential $\Gamma$ defined in \eqref{1.0}. Beyond the obvious regularity and symmetry properties, $\Gamma(x,\cdot)$ is harmonic in $\R^n-\{x\}$ and similarly for $\Gamma(\cdot,y)$, and \eqref{1.1} holds. In addition, if a function $u$ is defined and continuous in a neighborhood of $x\in\R^n$ then
\begin{align}\begin{split}\label{2.2}
&\lim_{\ep\to 0}\int_{\partial B_\ep(x)}\Gamma(x,z)u(z)\,dS_z =0,\\ &\lim_{\ep\to 0}\int_{\partial\eps B(x)}(\dnz\Gamma)(x,z)u(z)\,dS_z =u(x).
\end{split}\end{align}
Applying these facts together with the divergence theorem we arrive at the fundamental representation formula
\begin{align}\begin{split}\label{2.3}
u(x)=\int_{\bset}\Big[(\dnz \Gamma(x,z))u(z)-\Gamma(x,z)&\dn u(z)\Big]dS_z\\
&+\int_\set \Gamma(x,z)\Delta u(z)\,dz
\end{split}\end{align}
for the restrictions $u$ of $C^2$ functions on $\R^n$ to $\set$ and for $x\in\set$. The Green's functions $G_D$ and $G_N$ for $\set$ are then constructed by adding to $\Gamma$ functions $h_D$ and $h_N$ defined on $\set\times\ov\set$ which for fixed $x\in \set$ satisfy $\Delta h_D(x,\cdot)=0$ and $\Delta h_N(x,\cdot)=-|\set|^{-1}$ in $\set$ and with boundary conditions adjusted so as to cancel one or the other of the summands in the boundary integral in \eqref{2.3}. It will then follow that
\begin{equation}\label{2.4}
u(x)=\int_{\bset}\Big[(\dnz G_D(x,z))u(z)\Big]dS_z
+\int_\set G_D(x,z)\Delta u(z)\,dz
\end{equation}
and
\begin{equation}\label{2.5}
u(x)-\ov u=-\int_{\bset}\Big[G_N(x,z)\dn u(z)\Big]dS_z
+\int_\set G_N(x,z)\Delta u(z)\,dz.
\end{equation}
(Observe that \eqref{2.5} cannot hold for $u\equiv 1$ without the term $\ov u$ on the left.)
These constructions are made precise in the following theorem:

\begin{theorem} Let $\set$ be a $C^{k+1}$ domain in $\R^n$ in the sense of {\rm Definition~2.1} where $n\ge 3$ and $k\ge 1$. 
\smallskip 
                                                                                  
\setlength{\hangindent}{18pt}
\noindent {\rm (a)} There are functions $h_D$ and $h_N$ in $C^\infty(\set\times\set)$ such that the mappings $x\in\set\mapsto h_{D,N}(x,\cdot)$ and $y\in\set\mapsto h_{D,N}(\cdot,y)$ are continuous from $\set$ into $H^{k+1}(\set)$  and for fixed $x\in\set$ \medskip

\begin{equation} \label{2.6} \left\{ \begin{array}{l}
  \Delta_y h_D(x,y) = 0, \ y\in \set,\\
	\smallskip h_D(x,\cdot)=-\Gamma(x,\cdot)\ in\ L^2(\bset)\\
\end{array}\right.
 \end{equation}
 and
 \begin{equation} \label{2.7} \left\{ \begin{array}{l}
  \Delta_y h_N(x,y) = -|\set|^{-1},\ y \in\set,\\    
  
	\smallskip (\dn h_N)(x,\cdot)=-(\dn \Gamma)(x,\cdot)\ in\ L^2(\bset),\\
	
	\smallskip \ov{ h_N(x,\cdot)}=0.
	
\end{array}\right.
 \end{equation}
 
Also, $h_D$ and $h_N$ are symmetric: $h_{D,N}(x,y)=h_{D,N}(y,x)$ for $\nobreak{x,y\in\set}$.
\smallskip
 
\setlength{\hangindent}{18pt}
\noindent {\rm (b)} The functions $G_D \equiv \Gamma+h_D$ and $G_N(x,\cdot)\equiv\Gamma(x,\cdot)-\ov{\Gamma(x,\cdot)}+h_N(x,\cdot)$ satisfy the following:

\smallskip

{\rm($i$)} The representations \eqref{2.4} and \eqref{2.5} hold for the restrictions $u$ of $C^2$ functions on $\R^n$ to $\ov\set$.

\smallskip
{\rm($ii$)} $G_{D,N}(x,y)=G_{D,N}(y,x)$ for $(x,y)\in \set\times\set-\sd$.

\medskip\setlength{\hangindent}{24 pt}
{\rm($iii$)} Given $p\in [1,n/(n-2))$ there is a constant $C(\set,k,p)$ such that for all $x\in\set$, $$|G_{D,N}(x,\cdot)|_{L^p(\set)}\le C(\set,k,p).$$
\medskip

{\rm($iv$)} $G_D$ and $G_N$ are unique in the following sense: if $g\in H^2(\set)$ and if for some $x\in\set$ and for all $u$ as in \eqref{2.4} and \eqref{2.5}
\begin{equation}\label{2.8}
u(x)=\int_\bset (\dn g(y))u(y)\,dS_y +\int_\set g(y)\Delta u(y)\,dy,
\end{equation}
then $g=G_D(x,\cdot)$ as elements of $L^p$ for $p\in [1,n/(n-2))$. 
Similarly, if $g\in H^1(\set)$ with $\ov g=0$ and if for some $x\in\set$ 
\begin{equation}\label{2.85}
u(x)-\ov u=-\int_\bset g(y)\dn u(y)\,dS_y +\int_\set g(y)\Delta u(y)\,dy
\end{equation}
for all $u$ as above,
then $g=G_N(x,\cdot)$ as elements of $L^p$ for $p\in [1,n/(n-2))$.

\end{theorem}

\bigskip

\noindent The solutions $h_D$ and $h_N$ of \eqref{2.6} and \eqref{2.7} are understood in the usual weak sense, following the formalism in \cite{hoffbook} section 4 or \cite{hoff3d} section~2, (5)-(7) for example. Thus $h_{D}(x,\cdot)$ and $h_{N}(x,\cdot)$ are in $H^1(\set)$ and the weak forms of the Poisson equations in (a) hold for test functions in $H^1_0(\set)$ and $H^1(\set)$ respectively. Note that since $h_D$ and $h_N$ are continuous mappings from $\set$ into $H^{k+1}(\set)$ and since $k+1\ge 2$, their traces and the traces of their gradients are continuous mappings from $\set$ into $L^2(\bset)$.  

Various properties of $G_D$ and $G_N$ are easily derived from those of $\Gamma,\, h_D$ and $h_N$. For example, for $x\in\set$, $G_D(x,\cdot)$ is harmonic in $\set-\{x\}$ and $\Delta G_N(x,\cdot)=-1/|\set|$ pointwise in {\nobreak $\set-\{x\}$}. Also, $G_{D,N}(x,\cdot)\in H^{k+1}(\set-B_\ep(x))$ for all $\ep>0$, modulo equivalence classes, so that the traces of $G_D(x,\cdot)$ and $G_N(x,\cdot)$ and their gradients are in $L^2(\bset)$. In particular, $G_D(x,\cdot)=0$ and  $\dn G_N(x,\cdot)=0$ in $L^2(\bset)$.

Note that the bound in (b)($iii$) is independent of $d_\partial(x)$. This bound will serve as the essential starting point for the derivation of pointwise bounds in Theorem~4.1 for the derivatives of the Green's functions, bounds which will also be independent of distances to the boundary.
\bigskip

\noindent{\bf Proof of Theorem  3.1.} We begin with the constructions of the Green's functions following the exposition and notations in \cite{robert}.  Let $\eta$ be a smooth, nondecreasing function of $s\in \R$ such that $\eta(s)=0$ for $s\le 1/3$ and $\eta(s)=1$ for $s\ge 2/3$. Then for fixed $x\in \set$ define

$$f(x,y) =  \eta\Big(\frac{|x-y|}{d_\partial(x)}\Big)\Gamma(x,y).$$
Then $\Delta_yf(x,\cdot)\in C^\infty(\R^n)$ and the map $x\mapsto f(x,\cdot)$ is continuous from $\set$ into $H^j(\set)$ for all $j\ge 0$. We can therefore fix $x\in\set$ and let $w_D(x,\cdot)$ and $w_N(x,\cdot)$ be the solutions of the (formal) systems
\begin{equation} \label{3.1} 
 \left\{ \begin{array}{l}
  \Delta_y w_D(x,y) = \Delta_yf(x,y),\ y\in\set,\\
	\smallskip w_D(x,y)=0,\ y\in\bset\\
	\end{array}\right.
 \end{equation}
  and
\begin{equation} \label{3.2}  \left\{ \begin{array}{l}
  \Delta_y w_N(x,y) = \Delta_yf(x,y)-\ov {\Delta_yf(x,\cdot)},\ y\in\set,\\
	\smallskip D_{\nu_y} w_N(x,y)=0,\  y\in\bset,\\
	\smallskip \ov {w_N(x,\cdot)} =0.
\end{array}\right.
 \end{equation}
Observe  that the first two conditions in \eqref{3.2} could be inconsistent without the term $\ov{\Delta_y f(x,\cdot)}$. Note also that
$$\ov{\Delta_y f(x,\cdot)}=|\set|^{-1}\int_\set \Delta_y f=|\set|^{-1}\int_\bset(\dn \Gamma)(x,\cdot)=|\set|^{-1},$$
the last equality following by putting $u\equiv 1$ in \eqref{2.3}. The solutions $w_D$ and $w_N$ here are understood in the usual weak sense as described above. Thus $w_D(x,\cdot)\in H^1_0(\set)$ and $w_N(x,\cdot)\in H^1_{1^\perp}(\set)\equiv \{v\in H^1(\set ): \ov v=0\}$, and the weak forms of the Poisson equations above hold for test functions in $H^1_0$ and $H^1_{1^\perp}$ respectively. Also, Theorem 4.2 of \cite{hoffbook} shows that $w_D$ and $w_N$ are continuous mappings from $\set$ into $H^{k+1}(\set)$, and since $k+1\ge 2$, their traces and the traces of their gradients are continuous mappings from $\set$ into $L^2(\bset)$. All these statements are therefore true for the functions 
$$h_D(x,\cdot)\equiv w_D(x,\cdot)-f(x,\cdot)$$
 and 
 $$h_N(x,\cdot)\equiv w_N(x,\cdot)-f(x,\cdot) +\overline {f(x,\cdot)},$$ which are weak solutions of the systems \eqref{2.6} and \eqref{2.7}. In particular, $h_D(x,\cdot)$ is harmonic in $\set$ as is $h_N(x,\cdot)$ modulo a polynomial; these functions are therefore in $C^\infty(\set)$. 

Next we check that for fixed $y\in \set$, $h_{D}(x,y)$ and $h_{N}(x,y)$ are continuous functions of $x\in\set$. To see this for the Neumann case we take $u=h_N(\cdot,y)$ in the representation formula \eqref{2.3} (which is allowed because the set of restrictions to $\set$ of $C^2$ functions on $\R^n$ is dense in $H^2(\set)$ ). Substituting from \eqref{2.7} into the second and third terms on the right in \eqref{2.3} we see by inspection that these terms are continuous in $x$. For the first term on the right we have that $|(\dnz \Gamma)(x,z)|$ is bounded by a constant depending on $d_\partial(x)$ so that the difference between the values of this term at points $x_1$ and $x_2$ in $\set$ is bounded by
$$C|\nabla_z h_N(x_1,\cdot)-\nabla_zh_N(x_2,\cdot)|_{L^1(\bset)}\le C|h_N(x_1,\cdot)-h_N(x_2,\cdot)|_{W^{2,1}(\bset)},$$
which goes to zero as $x_2\to x_1$ because the map $x\mapsto h_N(x,\cdot)$ is continuous from $\set$ into $H^2(\set)$. The proof for $h_D$ is similar. Summarizing, we now have that $h_{D,N}\in C(\set\times\set)$ and that for $x\in\set$, $h_{D,N}(x,\cdot)\in C^\infty(\set)$ and depends continuously on $x\in \set$ as an element of $H^{k+1}(\set)$. 

The representations \eqref{2.4} and \eqref{2.5} for $G_D$ and $G_N$ now follow easily from \eqref{2.3}, \eqref{2.6} and \eqref{2.7}. We can apply these to prove the symmetry property (b)($ii$), first for the Neumann case. Let $\varphi$ be a smooth function on $\set$ with compact support in $\set$ and solve the system
 
 $$ \left\{ \begin{array}{l}
  \Delta u = \varphi\ {\rm in}\ \set \\
  	\dn u=0\ {\rm on}\ \bset\\
		\ov u=0
	\end{array}\right.
$$
for $u$, and similarly with $\varphi$ and $u$ replaced by $\psi$ and $v$. Then the function $G_N(x,y)\varphi(x)\psi(y)$ is continuous on $\set\times\set-\sd$, hence measurable on $\set\times\set$, and is locally integrable because $\Gamma$ is integrable on $\set\times\set$ and the $H^{k+1}$ norm of $h_N(x,\cdot)$ is bounded uniformly for $x$ in the support of $\varphi$. The Fubini and Tonelli theorems are therefore justified in the following computation in which the representation \eqref{2.5} is applied twice:

\begin{align*}
\int_{\set\times\set}&G_N(y,x)\varphi(x)\psi(y)\,d(x,y)=\int\Big(\int G_N(y,x)\varphi(x)\,dx\Big) \psi(y)\,dy \\
	=&\int u(y)\Delta v(y)dy=\int \Delta u(y)v(y) dy= \int \varphi(x)v(x)dx \\
	=&\int\varphi(x)\Big(\int G_N(x,y)\psi(y)\,dy\Big)\,dx=\int_{\set\times\set}G_N(x,y)\varphi(x)\psi(y)\,d(x,y).
\end{align*}
The symmetry statement in (b)($ii$) for $G_N$ now follows by letting $\varphi$ and $\psi$ tend to Dirac masses at distinct points
$x_0$ and $y_0$ in $\set$ and recalling that $h_N$ and therefore $G_N$ are continuous on $\set\times\set-\sd$. The proof for $G_D$ is similar. This symmetry together with the conclusion at the end of the previous paragraph then combine to complete the proofs of the regularity statements in part (a) of the theorem.

Next we prove the $L^p$ bound in (b){\it (iii)} for $p\in (1,n/(n-2))$. Let $p'\in (n/2,\infty)$ be its H\"older conjugate and let $\varphi\in L^{p'}(\set)$. Then for the Dirichlet case we solve the system 
$$ \left\{ \begin{array}{l}
  \Delta u = \varphi\ {\rm in}\ \set \\
  	u=0\ {\rm on}\ \bset
	\end{array}\right.
$$
for $u$, and similarly for the Neumann case with the boundary condition above replaced by the conditions $\dn u=0$ on $\bset$ and $\bar u=0$. Then $|u|_{W^{2,p'}}\le C|\varphi|_{L^{p'}}$ by standard elliptic regularity results (see \cite{gt} pg.~110 for example), and $\sup_\set|u|\le C|u|_{W^{2,p'}}$ since $p'>n/2$ (see \cite{adams} Theorem~4.12). Applying \eqref{2.4} for the Dirichlet case we then find that for $x\in\set$,
\begin{align*}\big|\int_\set |G_D(x,y)\varphi(y)\,dy| &=  \big|\int_\set |G_D(x,y)\Delta u(y)\,dy\big|\\  &= |u(x)|
	\le C|u|_{W^{2,p'}}\le C|\varphi|_{L^{p'}}.
\end{align*}
 This proves that $|G_D(x,\cdot)|_{L^p}\le C(\set,k,p)$. The argument for $G_N$ is similar.  
 
 Finally we prove the uniqueness statement in (b)($iv$) for the Neumann case. Letting $g$ be as in the hypothesis we subtract \eqref{2.5} from \eqref{2.85} to find that
$$\int_\bset(G_N(x,\cdot)-g)\dn u\,dS+\int_\set(g-G_N(x,\cdot))\Delta u=0$$
for restrictions to $\set$ of functions $u\in C^2(\R^n)$, hence for $u$ in $W^{2,p'}(\set)$. Letting $\varphi$ be an arbitrarily element of $L^{p'}(\set)$ we can choose $u\in W^{2,p'}(\set)$ to satisfy $\Delta u=\varphi$ in $\set$, $\dn u=0$ on $\bset$ and $\ov u=0$. The result then shows that $g=G_N(x,\cdot)$ as elements of $L^p(\set)$. The proof for the Dirichlet case is similar.

\medskip
\rightline\qedsymbol
\bigskip

\section{Pointwise Bounds} 
In this section we derive bounds for the Green's functions and their derivatives at points $(x, y)\in\set\times\set-\sd$, bounds which are independent of $d_\partial(x)$ and $d_\partial(y)$ and which detail the rates of blowup as $y\to x$. In particular, we show in Theorem~4.1(a) below that the derivatives of $G_D$ and $G_N$ satisfy the same bounds \eqref{1.1} satisfied by the fundamental solution $\Gamma$ up to an order determined by the regularity of $\bset$. 

To formulate these results we recall the  following imbedding theorem, which will be applied extensively throughout: Let $k_0$ be the positive integer 

\begin{equation} \label{2.9}k_0= \left\{ \begin{array}{l}
 n/2+1\ \text{ if $n$ is even}\\
	\smallskip (n+1)/2\ \text{if $n$ is odd;}\\
\end{array}\right.
 \end{equation}
 and let  $\sv\subset\R^n$ be a bounded open set in $\R^n$ with a Lipschitz boundary; that is, $\sv$ satisfies the conditions in Definition 2.1 but with the $\psi^z$ assumed only to be Lipschitz continuous. Then if $k=k_0+k_1$ with $k_1\ge 0$, elements of  $H^{k}(\sv)$ have specific Lebesgue representations in $C^{k_1,\lambda_0}(\set)$ where $\lambda_0=1/2$ if $n$ is even and $\lambda_0$ is an arbitrary number in $(0,1)$ if $n$ is odd; and
for $u\in H^{k}(\set)$,
\begin{equation}\label{2.10}
|u|_{C^{k_1,\lambda_0}(\set)}\le C(\sv,k,\lambda_0)|u|_{H^{k}(\set)}
\end{equation}
where the constant $C$ depends on an upper bound for the Lipschitz constants for the functions $\psi^z$ in Definition 2.1 and on $\lambda_0$ if $n$ is even. See \cite{adams}, Theorem 4.12, for the complete statement and related discussion.
\smallskip

Pointwise bounds for the Green's functions $G_D$ and $G_N$ are detailed in the following theorem:

\begin{theorem} Let $\set$ be a $C^{k+1}$ domain in $\R^n$ in the sense of {\rm Definition~2.1} where $n\ge 3$ and $k=k_0+k_1$ with $k_1\ge 0$, and let $\lambda_0$ be as above. \smallskip 
                                                                                  
\noindent There are constants $C(\set,k)$ and $C(\set,k,\lambda_0)$ such that for $|\alpha|,|\beta|\le k_1$ the derivatives $D^\alpha_xD_y^\beta G_N$ and $D^\alpha_xD_y^\beta h_N$exist in the sense of calculus on $\set\times\set-\sd$; these derivatives are H\"older-$\lambda_0$   
on $\set\times\set-\sd_\ep$ for every $\ep>0$ and therefore have continuous extensions to $\bar\set\times\bar\set-\sd_\ep$. 
The extensions satisfy 
\begin{equation}\label{2.11}
|D_x^\alpha D_y^\beta G_N(x,y)|,\,|D_x^\alpha D_y^\beta h_N(x,y)|\le C(\set,k)|x-y|^{(2-n)-(|\alpha|+|\beta|)}
\end{equation}
and
\begin{align}\begin{split}\label{2.115} 
\langle D_x^\alpha D_y^\beta G_N\rangle&_{(B_r(x)\cap\ov\set)\times( B_r(y)\cap\ov\set),\lambda_0}\\&\le C(\set,k,\lambda_0)|x-y|^{(2-n)-(|\alpha|+|\beta|)-\lambda_0}
\end{split}\end{align}
where $r(x,y)=\min\{|x-y|,1\}/C(\set,k)$.
Also, for $|\alpha|\le k_1$ and $|\beta|\le k$,  
\begin{equation}\label{2.12}
|D_x^\alpha D_y^\beta G_N(x,\cdot)|_{L^2(B_{r}(y)\cap\set)} \le C(\set,k)|x-y|^{(2-n)-(|\alpha|+|\beta|)}r^{n/2},
\end{equation}
for $(x,y)\in\ov\set\times\ov\set-\sd$.
\medskip

\noindent The same results hold for $G_D$ and $h_D$ but under the stronger hypothesis that $k_1\ge 1$. \end{theorem}
\medskip

The proof of Theorem 4.1 begins below, following the statement and proof of Corollary 4.2. Observe that the rates of blowup as $y\to x$ in the above bounds are exactly the same as in \eqref{1.1} for the fundamental solution $\Gamma$. The apparent discrepancy between the $C^{k_0+k_1+1}$ regularity assumed here for $\set$ and the $C^{k_1,\lambda_0}$ regularity proved for the Green's functions arises because the latter is derived from the $H^{k_0+k_1}$ bound in \eqref{2.12} by way of the imbedding result \eqref{2.10}.  On the other hand the $H^{k_0+k_1}$ regularity of $D_x^\alpha G_N(x,\cdot)$ in \eqref{2.12} is of importance in itself, and this entails a loss of only one order, arising in the construction of a certain cut-off function in Lemma~7.1 in the Appendix.

Bounds for the Green's functions which depend on distances to $\bset$ are also of interest and are easily derived from the maximum principle for harmonic functions, standard interior bounds (Theorem 7.3) and the bound in \eqref{2.11} for $\alpha=\beta=0$. We include these here both for completeness and for their usefulness in applications, including the application given in Section~6. Before giving the statement we recall that if $\set$ is a bounded open set in $\R^n$ then given $f\in C(\bset)$ there is a harmonic extension of $f$ to $\ov\set$ provided that there is a ``barrier function" at each point of $\bset$. This barrier condition is satisfied for example if $\set$ satisfies the conditions in Definition 2.1 with Lipschitz continuous mappings $\psi^y$; the barrier condition is much more general, however. See \cite{evans} pg. 367 or \cite{gt} Theorem 2.14 for detailed statements and related discussion.
\medskip
\begin{corollary} Let $\set$ be a bounded open set in $\R^n$, $n\ge 3$, and assume that there is a barrier at each point of $\bset$. Then the function $h_D$ defined in  \eqref{2.6} is in $C^\infty(\set\times\set)$ and satisfies
\begin{equation}\label{2.116}
|D_x^\alpha D_y^\beta h_D(x,y)|\le C[\db(x)^s\db(y)^{1-s}]^{(2-n)}\db(x)^{-|\alpha|}\db(y)^{-|\beta|}
\end{equation}
for all $(x,y)\in\set\times\set$, all $\alpha$ and $\beta$ and all $s\in [0,1]$, where $C$ depends only on $|\alpha|$ and $|\beta|$. The same bound holds for $h_N$ but under the stronger hypothesis that $\set$ is a $C^{k_0+1}$ domain and $C=C(\set,k_0+1,|\alpha|,|\beta|)$. 
\end{corollary}

\begin{proof} If $x,\,x_0\in\set$ with $|x-x_0|\le\db(x_0)/2$, then $\db(x)\ge\db(x_0)/2$. And since $h_D(x,\cdot)=-\Gamma(x,\cdot)$ on $\bset$,
\begin{equation}\label{2.13}|h_D(x,y)|\le C|x-y|^{2-n}\le C\db(x_0)^{2-n}
\end{equation}
for $y\in\bset$ and therefore for all $y\in\ov\set$ by the maximum principle. This bound holds in particular on $B_{\db(y_0)}(y_0)$ for a given fixed $y_0\in \set$ so that by standard interior estimates (see Theorem 7.3) 
$$|D_y^\beta h_D(x,y_0)|\le C(|\beta|)\db(x_0)^{2-n}\db(y_0)^{-|\beta|}$$
 for all $x\in B_{\db(x_0)/2}(x_0).$ Consequently by standard interior estimates again, applied to the harmonic function $D_y^\beta h_D(\cdot,y_0)$,
 $$|D_x^\alpha D_y^\beta h_D(x_0,y_0)|\le C(|\alpha|,|\beta|)\db(x_0)^{(2-n)-|\alpha|}\db(y_0)^{-|\beta|}.$$
The same result holds with $x_0$ and $y_0$ reversed and $\alpha$ and $\beta$ reversed, and the two bounds can then be combined to obtain \eqref{2.116}.

For $h_N$ we apply \eqref{2.11} in Theorem 4.1 to obtain that 
$$|h_N(x,y)|\le C(\set,k_0+1)|x-y|^{2-n}\le C\db(x_0)^{2-n},$$
 for $y\in\bset$, as in \eqref{2.13}. The rest of the proof is then the same as for $h_D$. 
 \end{proof}
 \medskip
 
\noindent Note that \eqref{2.116} does not hold for $G_{D,N}$ in general because otherwise it would hold for $\Gamma$, which in general it does not. We also note that geometric combinations of \eqref{2.11} and \eqref{2.116} can be taken and are often useful, as will be seen in the proof of Theorem~5.1 in the next section.\bigskip

\noindent{\bf Proof of Theorem 4.1.} The proofs of \eqref{2.11}, \eqref{2.115} and \eqref{2.12} are straightforward when either $x$ or $y$ is far from the boundary or when $x$ and $y$ are far apart. These cases are covered in Lemma 4.3 below. The remaining case that $x$ and $y$ are close and are close to the boundary is considerably more technical, involving approximations of the Green's functions by the Green's functions for the half-space determined by the tangent plane at $y'$ together with appropriately scaled interior-type estimates at the boundary and a scaled version of the Sobolev imbedding $H^1\mapsto L^{2n/(n-2)}$. Statements and proofs of the latter results for general elliptic equations are given in Theorems~7.2 and 7.4 in the Appendix and will be cited as needed, and all the definitions and notations in Section 2 will be in force in this section without further reference.

The proof of Theorem 4.1 requires three preparatory lemmas, the first of which gives bounds for the Green's functions for the two cases referred to above as well as a separation of points construction which will be applied here and in sections ~5 and 7.
 
\medskip

\begin{lemma} Let $\set$ be a $C^{k+1}$ domain in $\R^n$ in the sense of {\rm Definition~2.1} where $k\ge 1$ and $n\ge 3$.
\smallskip

\setlength{\hangindent}{18pt}
\noindent {\rm (a)}  Given a compact set $K\subset\set$ there is a constant $C(\set,k,K)$ such that if $x\in K$  then $|h_{D,N}(x,\cdot)|_{H^{k+1}(\set)}\le C(\set,k)$; and if $y\in K$ then  $|h_{D,N}(\cdot,y)|_{H^{k+1}(\set)}\le C(\set,k)$.
\smallskip

\setlength{\hangindent}{18pt}
\noindent {\rm (b)} There are constants $R_1(\set,k)\in (0,\min\{R_0,\kappa_0\})$ and $C(\set,k)$ such that the following separation of points property holds: Given distinct points $x$ and $y$ in $\ov\set$ let $s(x,y)= \min\{R_1,|x-y|\}$ and define sets $\sv_0(x,y),\,\sv_1(x,y)$ and $\sv_2(x,y)$ as follows: 
\begin{align*}
&\text{if $d_\partial(y)\ge s/12$ then $\sv_0=B_{s/14}(y),\, \sv_1=B_{s/16}(y),\,\sv_2=B_{s/18}(y)$};\\
&\text{if $d_\partial(y)< s/12$ then $\sv_0=\sw^{y'}_{s/3}(y),\, \sv_1=\sw^{y'}_{s/4}(y),\,\sv_2=\sw^{y'}_{s/6}(y)$}
\end{align*}
where $y'$ is the point on $\bset$ closest to $y$ and the sets $\sw^{y'}_t$ are as defined in \eqref{2.1}. Then in either case dist$(x,\sv_0(x,y))\ge C^{-1}|x-y|$ and $B_{s/18}(y)\subset\sv_2(x,y)$.
\medskip

\setlength{\hangindent}{18pt}
\noindent {\rm (c)} Given $r>0$ there is a constant $C(\set,k,r)$ such that if $x$ and $y$ are points of $\set$ with $|x-y|\ge r$ then 
$$|h_{D,N}(x,\cdot)|_{H^k(\sv_2(x,y))},\,|G_{D,N}(x,\cdot)|_{H^k(\sv_2(x,y))}\le C(\set,k,r^{-1}).$$
\end{lemma}
\medskip
\begin{proof} Part (a) follows directly from the statement in Theorem 3.1(a) that the maps $x\to h_{D,N}(x,\cdot)$ and $y\to h_{D,N}(\cdot,y)$ are continuous from $\set$ into $H^{k+1}(\set)$. 

The statements in (b) are obvious for the case that $\db(y)\ge s/12$. For the other case we first observe that if $w$ and $z$ are points of the cylinder  $B_t(0)\times (-t,0)$ in $\R^n$ and $z=(0,z_n)$, then $|w-z|\le \sqrt{2}t$. It follows that if $R_1$ is sufficiently small depending on the moduli of continuity of the $\nabla\psi^z$ in Definition 2.1 then for $t\le R_1$ the following holds: if $y\in\set$ with $\db(y)\le\kappa_0$ then $|w-z|\le 2t$ for points $w$ and $z$ in $\sw^{y'}_t$ with $z=(0,z_n)$. The statements in (b) for the case that $\db(y)<s/12$ then follow easily.

To prove (c) we begin with the bound in Theorem 3.1(b)(iii) with $p$ the average of one and $n/(n-2)$ so that $|G_{D,N}(x,\cdot)|_{L^p(\sv_0(x,y))}\le C(\set,k)$. Then since $G_D(x,\cdot)$ is harmonic in $\sv_0$, as is $G_N(x,\cdot)$ modulo a polynomial, Theorem 7.2 applies to show first by repeated application of part (c) through sets interpolating $\sv_0$ and 
$\sv_1$ that $|G_{D,N}(x,\cdot)|_{L^2(\sv_1)}\le C(\set,k,r^{-1})$, and then by repeated application of part (b) through sets interpolating $\sv_1$ and $\sv_2$ that $|G_{D,N}(x,\cdot)|_{H^k(\sv_2)}\le C(\set,k,r^{-1})$. The same bound holds for $\Gamma(x,\cdot)$ by inspection, hence for $h_{D,N}(x,\cdot)$.

\end{proof}

In the following lemma we define the reflection $\ov y$ of a point $y$ near $\bset$ and we give an approximate Pythagorean theorem for the triangle whose vertices are $y'$, a nearby point $z\in\bset$ and either $y$ or $\ov y$: 
\medskip

\begin{lemma} Let $\set$ be a $C^{k+1}$ domain in $\R^n$ in the sense of {\rm Definition~2.1} where $k\ge 1$ and $n\ge 2$.
\smallskip

\setlength{\hangindent}{18pt}
\noindent {\rm (a)} If $y_0\in\set$ with $\db(y_0)\le\kappa_0$ then there is a positive $\ep$ depending on $d_\partial(y_0)$ such that for $y\in B_\ep(y_0)$, say $y=S^{y_0'}(\tilde w,w_n)$, the reflected point $
\ov y\equiv S^{y_0'}(\tilde w,-w_n)$ is in $S^{y'}(\cy_{R_0})-\ov\set$. 
\smallskip

\setlength{\hangindent}{18pt}
\noindent {\rm (b)} There is constant $R_2\in (0,R_1]$ depending on $\set$ and $k$ and a constant $C(\set,k)$ such that if $y_0$ and $\ov y_0 $ are as in {\rm (a)} and $z\in \bset\cap S^{y_0'}(\cy_{R_2})$, say $z=S^{y_0'}((\tilde w,\psi^{y_0'}(\tilde w))$, then
$$C^{-1}\big(|\tilde w|^2+\db(y_0)^2\big)\le |z-y_0|^2,\,|z-\ov y_0|^2\le  C\big(|\tilde w|^2+\db(y_0)^2\big).$$
\end{lemma}
\medskip

\begin{proof} Since $\db(y_0)\le \kappa_0<R_0$, $y_0$ is in $\cy^{y_0'}_{R_0}$, hence so is $y$ if $\ep$ is small. In particular, if $T^{y_0'}(y)=(\tilde w, w_n)$ then $|\tilde w|<R_0$ and $-R_0<w_n<R_0$. Also, since $T^{y_0'}(y_0)=(0, -\db(y_0))$,
$$|\psi^{y_0'}(\tilde w)|\le |\psi^{y_0'}(0)|+C|y-y_0|\le C\ep$$
and 
$$w_n\le (T^{y_0'}(y_0))_n+C|y-y_0| \le -\db(y_0) +C\ep.$$
Thus $$R_0>-w_n>\db(y_0)+\psi^{y_0'}(\tilde w)-C\ep>\psi^{y_0'}(\tilde w)$$
if $\ep$ is small. This proves that $\ov y\in S^{y_0'}(\cy_{R_0})-\ov\set=\cy^{y_0'}_{R_0}-\ov\set$. 

The result in (b) follows because 
$$|z-y_0|,\,|z-\ov y_0|=|(\tilde w,\psi(w))-(0,\mp\db(y_0)|$$
 and $\psi^{y_0'}(0)$ and $\nabla\psi^{y_0'}(0)$ are zero.
 
 \end{proof}
\medskip

\noindent Note that the definition $\ov y=S^{y_0'}(\tilde w,-w_n)$ in part (a) above shows that $\ov y$ depends on $y_0'$ as well as on $y$. We will therefore refer to $\ov y$ as {\it the reflection of $y$ relative to $y_0'$} when required for clarity.
\medskip

In the following lemma we give representations for $h_D(x,y)$ and $G_N(x,y)$ for $y$ in a neighborhood of a point $y_0$ near $\bset$ in terms of $\Gamma(y,\cdot)-\Gamma(\ov y,\cdot)$ and $\Gamma(y,\cdot)+\Gamma(\ov y,\cdot)$, these being the Dirichlet and Neumann Green's functions respectively for the half space determined by the plane tangent to $\bset$ at $y_0'$. 

\begin{lemma} Let $\set$ be a $C^2$ domain in $\R^n$ as in {\rm Definition 2.1} and let $R_2,\,y_0\in \set$ and 
$y\in B_\ep(y_0)$ be as in { \rm Lemma 4.4} and $x\in\set$.
\smallskip

\setlength{\hangindent}{18pt}\noindent {\rm (a)} If $\ov y$ is the reflection of $y$ relative to $y_0'$ then\smallskip
\begin{align} \begin{split}\label{2.131}
h_D(x,y)=-\int_\bset&\big(\dnz [\Gamma(y,z)-\Gamma(\ov y,z)]\big)\Gamma(x,z)dS_z\\
	&-\int_\bset[\Gamma(y,z)-\Gamma(\ov y,z)]\dnz h_D(x,z)dS_z;
\end{split}\end{align} 
and if $x\not=y$ then
\begin{align}\begin{split}\label{2.132}
\tilde G_N(x,y)=\Gamma(x,y)&+\Gamma(x,\ov y)-\ov{\Gamma(\cdot,y)} + \ov {\Gamma(\cdot,\ov y)}\\
	&+\int_\bset\big(\dnz [\Gamma(y,z)+\Gamma(\ov y,z)]\big)\tilde G_N(x,z)dS_z
\end{split}\end{align}
where $\tilde G_N=\Gamma+h_N$. 
\smallskip

\setlength{\hangindent}{18pt}
\noindent {\rm (b)} If $\set$ is a $C^{k+1}$ domain in $\R^n$ with $k\ge 1$ then for $z\in \bset\cap\ov{\sw^{y_0'}_{R_2}}$,
\begin{align}\begin{split}\label{2.133}&|\Gamma(y_0,z)-\Gamma(\ov y_0,z)|\le C\db(y_0)|z-y_0|^{2-n},\\
& |\dnz(\Gamma(y_0,z)+\Gamma(\ov y_0,z))|\le C|z-y_0|^{2-n}
\end{split}\end{align}
and for unit vectors $\sigma$, 
\begin{align}\begin{split}\label{2.134}|\sigma\cdot \nabla_z& \big(\Gamma(y_0,z)-\Gamma(\ov y_0,z)\big)|\\
&\le C\db(y_0)\begin{cases} |z-y_0|^{1-n}\ {\rm if}\ \sigma\perp\nu(y_0)\\ |z-y_0|^{-n}\ {\rm if}\ \sigma\parallel\nu(y_0). \end{cases}
\end{split}\end{align}

\end{lemma}
 
\begin{proof} 
To prove the representation in (a) for $h_D$ we note first that since $h_D(x,z)$ and $\Gamma(\ov y,z)$ are harmonic functions of $z\in\set$,
$$\int_\bset\big[\big(\dnz \Gamma(\ov y,z)\big)h_D(x,z)-\Gamma(\ov y,z)\dnz  h_D(x,z)\big]dS_z=0;$$
and from \eqref{2.3} 
$$h_D(x,y)=\int_\bset\big[(\dnz \Gamma(y,z))h_D(x,z)-\Gamma(y,z)\dnz h_D(x,z)\big]dS_z.$$
The representation for $h_D$ then follows by putting $h_D(x,\cdot)=-\Gamma(x,\cdot)$ on $\bset$ and subtracting.

To derive the representation for $\tilde G_N$ we recall from \eqref{2.7} that \break$\nobreak{\Delta_zh_N(x,z)=-|\set|^{-1}}$. Multiplying by  $\Gamma(\ov y,z)$, which is harmonic in $z$ for $z\in \set$, and integrating  we therefore get
$$\int_\bset\Big( \Gamma(\ov y,z)D_{\nu_z}h_N(x,z)-(D_{\nu_z}\Gamma(\ov y,z))h_N(x,z)\Big)dS_z=-\ov{\Gamma(\ov y,\cdot)}.$$
Subtracting this from the representation \eqref{2.4} with $u=h_N(x,\cdot)$ we then obtain
\begin{align}\begin{split}\label{2.14}
h_N(x,y)=\int_\bset\Big[(D_{\nu_z}[\Gamma(y,z)+\Gamma(\ov y,z)])h_N(x,z)&\\ -(\Gamma(y,z)+\Gamma(\ov y,z))D_{\nu_z}&h_N(x,z)\Big]dS_z\\
-\ov{\Gamma(\cdot, y)} +\ov{\Gamma(\cdot, \ov y)}.&
\end{split}
\end{align}
Next we derive a corresponding representation for $\Gamma(x,y)+\Gamma(x,\ov y)$. First, since $\Gamma(x,\ov y)$ is harmonic in $x\in\set$, we get from \eqref{2.4} that
$$\Gamma(x,\ov y)=\int_\bset\big[(D_{\nu_z}\Gamma(x,z))\Gamma(z,\ov y)-\Gamma(x,z)D_{\nu_z}\Gamma(\ov y,z)\big]dS_z.$$
Also, since $\Gamma(x,y)$ is a harmonic function of $x\in\bset-B_\delta(y)$ for small $\delta$, we can apply the representation \eqref{2.3} with $\set$ replaced by $\set-B_\delta(y)$ then let $\delta\to 0$ and apply \eqref{2.2} to get that
$$\int_\bset\big[(D_{\nu_z}\Gamma(x,z))\Gamma(y,z)-\Gamma(x,z)D_{\nu_z}\Gamma(y,z)\big]dS_z=0.$$
Adding the last two equations we obtain
\begin{align}\begin{split}\Gamma(x,\ov y)=\int_\bset\big[(D_{\nu_z}\Gamma(x,z))\big(\Gamma(y,z)+&\Gamma(\ov y,z)\big)\\
	-\Gamma(x,z)D_{\nu_z}&\big(\Gamma(y,z)+\Gamma(\ov y,z)\big)\big]dS_z.
\end{split}\end{align}
We subtract this from \eqref{2.14}, recalling that $D_{\nu_z}h_N(x,z)=-D_{\nu_z}\Gamma(x,z)$ for $z\in\bset$ to obtain finally that

\begin{align*} h_N(x,y)=\Gamma(x,\ov y)+\int_\bset&\big(D_{\nu_z}\big[\Gamma(y,z)+\Gamma(\ov y,z)\big]\big)\big(\Gamma(x,z) +h_N(x,z)\big)dS_z\\ &- \ov{\Gamma(\cdot,y)}+\ov{\Gamma(\cdot,\ov y)} .
\end{align*}
The representation in (a) for $\tilde G_N$ then follows by adding $\Gamma(x,y)$ to both sides.

To prove the bounds in (b) we recall that   $T^{y_0'}(y_0)=(0,-\db(y_0))$, $T^{y_0'}(\ov { y_0})=(0,\db(y_0))$ and $T^{y_0'}(z)=(\tilde w,\psi^{y_0'}(\tilde w))$ for some $\tilde w$. Then since $T^{y_0'}$ is an isometry,
$$\big| |y_0-z|^2-|\ov y_0-z|^2\big| \le C\db(y)|\tilde w|^2.$$ 
The required bounds then follow by elementary computations together with the results in Lemma 4.4(b) and the fact tht $\ov{y_0}-y_0=2\db(y_0)\nu(y_0).$

\end{proof}

\medskip

We now complete the proof of Theorem 4.1, fixing $k=k_0+k_1$ as in the hypothesis.

\medskip

\noindent{\bf Step 1: Proof of \eqref{2.11} with $\alpha=\beta=0$.} We will prove that there is a constant $C(\set,k)$ such that 
\begin{equation}\label{2.15}
|h_{D,N}(x,y)|,\,|G_{D,N}(x,y)| \le C(\set,k)|x-y|^{2-n}
\end{equation}
for all $(x,y)\in \set\times\set-\sd$.  This bound holds by Lemma~4.3(a) if either $\db(x)\ge\kappa_0$ or $\db(y)\ge \kappa_0$ since $H^{k+1}(\set)$ is continuously included in $C^{k_1}(\set)$ and $C$ may depend on the diameter of $\set$. We therefore consider only points $x$ and $y$ in $\set$ both of which are within $\kappa_0$ of $\bset$. We also note that if $|x-y|\ge R$ for some $R>0$ then Lemma~4.3(c) applies to show that the terms on the left in \eqref{2.15} are bounded by $C(\set,k,R^{-1})$.

We begin with the Neumann case, this being somewhat less technical.  
Fix $x\in\set$ with $\db(x)<\kappa_0$, temporarily fix an $R\in(0, R_2/4)$ and define
$$M_R=\sup\{|x-y|^{n-2}|\tilde G_N(x,y)| : y\in \ov{\set\cap B_R(x)}\ {\rm and}\ \db(y)\le \kappa_0\} $$
where $\tilde G_N=\Gamma+h_N$. Then $M_R$ is nonempty and the sup is finite by properties of 
$\Gamma$ and $h_D(x,\cdot)$. To bound $M_R$ we  fix a point $y_0$ in the defining set and may assume without loss of generality that $y_0\in\set\cap B_R(x)$ because $\Gamma(x,\cdot)$ and $h_N(x,\cdot)$ are continuous on $\ov\set-\{x\}$. 

Now consider the representation \eqref{2.132} for $\tilde G_N(x,y_0)$. The first four terms on the right are easily seen to be bounded by $C(\set,k)|x-y_0|^{2-n}$ and we claim that in the integral on the right the contribution from points $z$ in the set  $\bset\cap B_R(x)^c$ is bounded by $C(\set,k,R^{-1})$. To see this note that if $|z-x|\ge R$ then $ |\tilde G_N(x,z)|\le C(\set,k,R^{-1})$ and if also $R<|x-z|\le R_2/2$ then $|y_0-z|\le R_2 $ and in this case the second bound in Lemma~4.5(b) applies to show that 
\begin{equation}\label{2.16}
|\dnz [\Gamma(y_0,z)+\Gamma(\ov y_0,z)]|\le C(\set,k)|z-y_0|^{2-n}.
\end{equation}
If $|x-z|\ge R_2/2$ then $|y_0-z|\ge R_2/4$ and the same bound therefore holds. Thus the integrand in \eqref{2.132} is bounded in $\bset\cap B_R(x)^c$ by $C(\set,k,R^{-1})|z-y_0|^{2-n}$, whose integral is bounded by $C(\set,k,R^{-1})$. 

Next we consider the contribution to the integral in \eqref{2.132} from points $z\in \bset\cap B_R(x)$. For such $z$ we have that $|\tilde G_N(x,z)|\le M_R|z-x|^{2-n}$ and since $|z-y_0|\le R/2$ here, \eqref{2.16} again holds. Assembling these bounds, we conclude that
\begin{align*}|\tilde G_N(x,y_0)|\le &C(\set,k)|x-y_0|^{2-n} + C(\set,k,R^{-1})\\
	&+C(\set,k)M_R\int_{\bset\cap B_R(x)}|z-x|^{2-n}|z-y_0|^{2-n}\,dS_z.
\end{align*}
We claim that the integral $I$ here is bounded by $C(\set,k)R|x-y_0|^{2-n}$. To see this we apply the triangle inequality to 
bound $|x-y_0|^{n-2}I$ by the integral over $\bset\cap B_R(x)$ of $|z-x|^{2-n} + |z-y_0|^{2-n}$ modulo a multiplicative constant. 
The first of these can be computed explicitly in local polar coordinates and is seen to be bounded by $C(\set, k)R$. For the second we simply note that if $|z-x|\le R$ then $|z-y_0|\le 2R$, so that its integral satisfies the same bound. Multiplying both sides above by $|x-y_0|^{n-2}$ we thus obtain
\begin{equation}\label{2.161}
|x-y_0|^{n-2}|\tilde G_N(x,y_0)|\le C(\set,k,R^{-1})+C(\set,k)RM_R,
\end{equation}
and taking the sup over $y_0$ in the set defining $M_R$ that 
\begin{equation}\label{2.17} 
M_R\le C(\set,k,R^{-1})+C(\set,k)RM_R.
\end{equation}
We choose $R(\set,k)=C(\set,k)/2$ and conclude finally that $M_R\le C(\set,k)$ for a new constant $C(\set,k)$. And since $|G_N(x,y)-\tilde G_N(x,y)|\le C(\set,k)$ by Theorem 3.1(b), the bound in \eqref{2.15} holds for $G_N$, hence for $h_N$. 

The proof of \eqref{2.15} for $G_D$ and $h_D$ is somewhat more technical. Again fix $x\in\set$ with $\db(x)<\kappa_0$, temporarily fix $R\in (0,R_2/4]$ and define (a new)
\begin{align}\begin{split}
\label{2.175}
M_R=\sup\{|x-y|^{n-2}|h_D(x,y)&| : y\in \ov{\set\cap B_R(x)}\\
& \text{and}\ 0\le\db(y)\le\db(x)\}.
\end{split}\end{align}
The above set is clearly nonempty and the sup is finite by Theorem~3.1(a). To bound $M_R$ we again fix $y_0$ in the above set, assuming without loss of generality that $y_0\in\set$. Then by \eqref{2.131}
\begin{align}\begin{split}\label{2.18}
-h_D(x,y_0)=\int_\bset(\dnz&\tilde\Gamma(z))\Gamma(x,z)dS_z\\
	&+\int_\bset\tilde\Gamma(z)\dn_z h_D(x,z)dS_z
\end{split}\end{align}
where $\tilde\Gamma(z)=\Gamma(y_0,z)-\Gamma(\ov y_0,z)$. We claim that the first integral on the right here is bounded by $C(\set,k)|x-y_0|^{2-n}$. To see this we let $\delta>0$ be small and integrate by parts over the set $\set-B_\delta(x)-B_\delta(y_0)$ then let $\delta\to 0$ and apply \eqref{2.2} to obtain 
 $$\int_\bset(\dnz\tilde\Gamma(z))\Gamma(x,z)dS_z=\Gamma(y_0,x)-\Gamma(\ov y_0,x)-\int_\bset \tilde\Gamma(z)\dnz\Gamma(x,z)dS_z.$$
 The first two terms on the right are bounded by $C(\set,k)|x-y_0|^{2-n}$ because $|x-\ov y_0|\ge C^{-1}|x-y_0|$ if $R$ is small depending on $\set$. Also, by Lemma 4.5(b) the integrand on the right is bounded by 
 $$C\db(y_0)|z-y_0|^{2-n}|z-x|^{1-n}\le C|z-y_0|^{2-n}|z-x|^{2-n}$$
  because $\db(y_0)\le \db(x)\le |z-x|$. The resulting integral is then bounded by $C(\set,k)|x-y_0|^{2-n}$ as in the discussion above preceding \eqref{2.161}. 
  
Next we bound the second integral on the right in \eqref{2.18}. If $z\in\bset-B_{R/2}(x)$ then $|\nabla_zh_D(x,z)|\le C(\set,k,R^{-1})$ by Lemma 4.3(c), \eqref{2.10} and our assumption that $k_1\ge 1$. Therefore since $\tilde\Gamma$ is integrable on $\bset$, the contribution to this integral from points in $\bset-B_{R/2}(x)$ is bounded by $C(\set,k,R^{-1})$. Assembling these bounds we then have from \eqref{2.18} that
\begin{align}\begin{split}\label{2.19} 
|h_D(x,y_0)|\le C(\set,k)|x-&y_0|^{2-n}+C(\set,k,R^{-1}) \\+\Big|&\int_{\bset\cap B_{R/2}(x)} \tilde\Gamma(z)\dn_z h_D(x,z)dS_z\Big|.
\end{split}
\end{align}

The next step is to bound the term $\dnz h_D(x,z)$ in the integrand on the right here in terms of $M_R$. To do this we fix $z\in \bset\cap B_{R/2}(x)$ and first bound the $L^2$ norm of $h_D(x,\cdot)$ in a small neighborhood of $z$. Define $t=\min\{|z-x|/2,\,\db(x)\}$ and check that if $w\in\set\cap B_t(z)$ then $w\in \set\cap B_R(x)$ and $\db(w)\le \db(x)$. The point $w$ is thus in the set defining $M_R$ and therefore $|h_D(x,w)|\le M_R|w-x|^{2-n}\le CM_R|z-x|^{2-n}$. Then since $\sw^z_{t/2}\subset \set\cap B_t(z)$, 
$$|h_D(x,\cdot)|_{L^2(\sw^z_{t/2})}\le C(\set,k)M_R|x-z|^{2-n}t^{n/2}.$$
Next, since $h_D(x,\cdot)$ is harmonic in $\set$ we can apply Theorem 7.2(b) to find that for $j\le k=k_0+k_1$,
$$\sum_{|\alpha|=j}|D^\alpha h_D(x,\cdot)|_{L^2(\sw^z_{t/4})}\le C(\set,k)M_R|z-x|^{2-n}t^{n/2-j}$$
and then by the scaled imbedding result $H^{k_0+1}\mapsto C^{1,\lambda_0}$ in Theorem~7.4(a) that
\begin{align*}
\text{sup}_{w\in \sw^z_{t/4}}|\nabla h_D(x,w)|\le C(\set,k)\}&\sum_{j=0}^{k_0}t^{j-n/2}|D^{j+1}h_D(x,\cdot)|_{L^2(\sw^z_{t/4})}\\
	\le C(\set,k)M_R&t^{-1}|z-x|^{2-n}\\
	\le C(\set,k)M_R&(\db(x))^{-1}|z-x|^{2-n}
\end{align*}
since $t\ge C^{-1}\db(x)$. In particular, $\nabla_zh_D(x,z)$ satisfies this same bound for all $z$ in the integrand in \eqref{2.19}. Combining this with the bound in Lemma 4.5(b) for $|\tilde\Gamma(z)|$, we conclude that the  
integrand in \eqref{2.19} is bounded by $C(\set,k)RM_R|z-x|^{2-n}|z-y_0|^{2-n}$ and therefore that the integral is bounded by $C(\set,k)RM_R|x-y_0|^{2-n}$, as in the discussion preceding \eqref{2.161}. Applying this bound in \eqref{2.19} and then multiplying both sides by $|x-y_0|^{n-2}$ we conclude that
\begin{equation}\label{2.195}|x-y_0|^{n-2}| h_D(x,y_0)|\le C(\set,k,R^{-1})+C(\set,k)RM_R.
\end{equation}
The proof now proceedes just as for the Neumann case starting from \eqref{2.161}, thus proving that \eqref{2.15} holds for all $x$ and $y$ in $\set$ with $0<\db(y)\le\db(x)<\kappa_0$. The condition that $\db(y)\le\db(x)$ is eliminated by symmetry, and cases in which either $\db(x)$ or $\db(y)$ exceeds $\kappa_0$ were included in the discussion at the beginning of the proof of Step 1, which is now complete. 
\bigskip

\noindent{\bf Step 2: Bounds for derivatives.} The next step is to bound $|D_y^\alpha G_{D,N}|$ in $\set\times\set-\sd$. As at the beginning of Step 1 we may restrict consideration to points $x$ and $y$ in $\set$ which are close and which are close to $\bset$. For such $x$ and $y$ let $\sv_i=\sv_i(x,y)$ as in Lemma~4.3(b) and recall that $\sv_0\supset\sv_1\supset\sv_2,$ dist$(x,\sv_0)\ge C^{-1}|x-y|,\, B_{s(x,y)/18}(y)\subset \sv_2$ and 
$|\sv_i|\le C(\set,k)|x-y|^{n}$. We have that $G_D(x,\cdot)$ is harmonic in $\sv_0$, as is $G_N(x,\cdot)$ modulo a polynomial, and $|G_{D,N}(x,\cdot)|_{L^2(\sv_0)}\le C|x-y|^{(2-n)+n/2}$ by \eqref{2.15}. Theorem~7.2(b) therefore applies to show that \begin{equation}\label{2.20}
|D^\gamma_yG_{D,N}(x,\cdot)|_{L^2(\sv_2)}\le C|x-y|^{(2-n)+n/2-|\gamma|}
\end{equation}
for $|\gamma|\le k=k_0+k_1$. Then by \eqref{a.61} in Theorem 7.4(a), 
\begin{align*}{\rm sup}_{\sv_2}|D_y^\alpha G_{D,N}(x,\cdot)|&\le C\sum_{i=0}^{k_0}\sum_{|\beta|=i} |x-y|^{i-n/2}|D_y^{\alpha+\beta} G_{D,N}(x,\cdot)|_{L^2(\sv_2)}\\
&\le C(\set,k)|x-y|^{(2-n)-|\alpha|}
\end{align*}
for $|\alpha|\le k_1$; and by \eqref{a.62}, 
 $$\langle D^\alpha_yG_{D,N}(x,\cdot)\rangle_{C^{0,\lambda_0}(\sv_2)}\le C(\set,k,\lambda_0)|x-y|^{2-n-|\alpha|-\lambda_0}.$$
Summarizing and applying symmetry, we thus have that for $(x,y)\in \set\times\set-\sd$ and for $|\alpha|\le k_1$, 
\begin{equation}\label{2.21}
|D_y^\alpha G_{D,N}(x,y)|,\, |D_x^\alpha G_{D,N}(x,y)|\le C(\set,k)|x-y|^{(2-n)-|\alpha|}
\end{equation}
and 

\begin{align}\begin{split}\label{2.22}
\langle D_y^\alpha G_{D,N}(x,\cdot)\rangle_{C^{0,\lambda_0}(\sv_2(x,y))},\, \langle &D_x^\alpha G_{D,N}(\cdot,y)
\rangle_{C^{0,\lambda_0}(\sv_2(y,x))} \\
&\le C(\set,k,\lambda_0)|x-y|^{2-n-|\alpha|-\lambda_0}.
\end{split} \end{align}

Next we consider mixed derivatives. Again fix $x$ and $y$ and $\sv_i(x,y)$ as above. Then by \eqref{2.21}
$$|D_x^\alpha G_{D,N}(x,\cdot)|_{L^2(\sv_0(x,y))}\le C(\set,k)|x-y|^{(2-n)+n/2-|\alpha|}$$
for $|\alpha|\le k_1$.  Also, an easy argument based on test functions shows that if $\alpha\not= 0$ then $D_x^\alpha G_{D,N}(x,\cdot)$ is harmonic in $\set-\{x\}$. We can therefore apply Theorem 7.2 again, exactly as in the derivation of \eqref{2.20} above to obtain that
\begin{equation}\label{2.23}
|D^\gamma_yD_x^\alpha G_{D,N}(x,\cdot)|_{L^2(\sv_2(x,y))}\le C|x-y|^{(2-n)+n/2-(|\alpha|+|\gamma|)}
\end{equation}
for $|\gamma|\le k$, and then by Theorem 7.4 that
\begin{equation}\label{2.24}
|D_y^\beta D_x^\alpha G_{D,N}(x,y)|\le C(\set,k)|x-y|^{(2-n)-(|\alpha|+\beta|)}
\end{equation}
and 
\begin{align}\begin{split}\label{2.25}
\langle D_y^\beta D_x^\alpha G_{D,N}(x,\cdot)\rangle_{C^{0,\lambda_0}(\sv_2(x,y))},\, \langle &D_y^\beta D_x^\alpha G_{D,N}(\cdot,y)
\rangle_{C^{0,\lambda_0}(\sv_2(y,x))} \\
\le C(\set,k,\lambda_0)&|x-y|^{2-n-(|\alpha|+|\beta|)-\lambda_0}
\end{split} 
\end{align}
for $|\beta|\le k_1$.
This proves the bounds in \eqref{2.11} and \eqref{2.115} for points $(x,y)\in\set\times\set-\sd$.
\medskip

\noindent{\bf Step 3: Extensions to the Boundary.} Next we extend the bounds \eqref{2.24} and \eqref{2.25} to points on $\bset$. Let $|\alpha|,\,|\beta|\le k_1$, fix $(x,y)\in \break\ov\set\times\ov\set-\sd$ and suppose that $(x_i,y_i),\,i=1,2,$ are points of $\set\times\set-\sd$ which are within $s(x,y)/C$ of $(x,y)$, where $s(x,y)$ is as in Lemma~4.3(b). Then $|x_i-y_i|\ge C^{-1}|x-y|$ if $C$ is large so that \eqref{2.24} and \eqref{2.25} hold with $x$ and $y$ replaced by $x_i$ and $y_i$ on the left but not on the right. Also, if $C$ is sufficiently large, 
$$|y_2-y_1|\le s(x,y)/C\le s(x_2,y_2)/18$$ 
so that by Lemma 4.3(b), $y_1\in \sv_2(x_2,y_2)$ and therefore by \eqref{2.25} with $(x,y)$ replaced by $(x_2,y_2)$ on the left,
  $$\Big|D^\beta_yD_x^\alpha G_{D,N}(x_2,\cdot)\big|_{y_1}^{y_2}\Big|\le C|x-y|^{2-n-(|\alpha|+|\beta|)-\lambda_0}|y_2-y_1|^{\lambda_0};$$
and similarly
$$\Big|D^\beta_yD_x^\alpha G_{D,N}(\cdot, y_1)\big|_{x_1}^{x_2}\Big|\le C|x-y|^{2-n-(|\alpha|+|\beta|)-\lambda_0}|x_2-x_1|^{\lambda_0}.$$ 
Triangulating, we conclude that for $C=C(\set,k, \lambda_0)$ sufficiently large, $\Big|D^\beta_yD_x^\alpha G_{D,N}\big|_{(x_1,y_1)}^{(x_2,y_2)}\Big|$ is bounded by the sum of the right sides above and therefore that $D^\beta_yD_x^\alpha G_{D,N}$ extends H\"older-continuously to the closure $(\ov{B_{s/C}(x)\cap\set)}\times \ov{(B_{s/C}(y)\cap\set)}$ with constant 
$$C(\set,k,\lambda_0)|x-y|^{2-n-(|\alpha|+|\beta|)-\lambda_0}.$$ The more convenient formulation \eqref{2.115} then follows because $s/C=\min\{|x-y|,1\}/C$ for a new constant $C$. The statements in \eqref{2.11}, \eqref{2.115} and \eqref{2.12} then follow from this H\"older continuity and \eqref{2.23}. This completes the proof of Theorem~4.1.
\medskip

\rightline\qedsymbol

\section{The Cancellation Property}
 \bigskip
 
 In this section we formulate and prove the cancellation property described in the introduction, this being the main result of the paper. Definitions and notations from section 2 will be in force and will be applied throughout. We  recall in particular the sets $\scy^z_R$ and $\sw^z_R$ in Definition~2.1 for points $z\in\bset$ and the constants $R_2\le R_1\le R_0$ and $\kappa_0\le R_0$ defined in Definition~2.1, Lemma 4.3(b) and Lemma 4.4(b).
 
Before formulating the theorem we introduce local tangential and normal vector fields $\tau^z_i$ and $N^z$ as follows. Let $z\in\bset$ and let $T^{z}$ be the rigid motion described in Definition~2.1, say $T^{z}$ is multiplication by an orthogonal matrix $Q^{z}$ followed by a translation. Then for $y\in \scy^{z}_{R_0}$ and $w=(\tilde w,w_n)=T^{z}(y)$ define vector fields $\tau^{z}_i(y)=(Q^{z})^{tr}(e_i+\psi_{w_i}(\tilde w)e_n)$ for $i\le n-1$, where $e_i$ the standard basis vector whose $j$-th component is the Kronecker delta $\delta^{i,j}$; and $N^{z}(y)=(Q^{z})^{tr}(-\psi^{z}_{w_1},\ldots, -\psi^{z_0}_{w_{n-1}},1)$. Thus if $y\in\bset\cap\scy^{z}_{R_0}$ then the unit outward normal vector at $y$ is $\nu(y)=N^{z}(y)/|N^{z}(y)|$ and the $\tau_i^{z}(y)$ span the tangent space to $\bset$ at $y$. 

In part (a) of the following theorem we state the cancellation property for directional derivatives of the Green's functions determined locally by the $\tau^z_i$. Parts (b) and (c) are reformulations which may be useful in particular applications. The proof of part (a) occupies the greater part of this section and is similar to but considerably more technical than that of Theorem 4.1, relying again on local approximation by Green's functions for half-spaces determined by planes tangent to $\bset$. 	
\bigskip

\begin{theorem} Let $\set$ be a $C^{k+1}$ domain in $\R^n$ in the sense of {\rm Definition~2.1} where $n\ge 3$ and $k=k_0+k_1$ with $k_1\ge 1$, and let $\tau^z_i$ and $N^z$ be as above for $z\in \bset$.  Then there is a constant $C(\set,k)$ and a positive number $R_3(\set,k)\le R_2$ such that the following hold:
\smallskip

\setlength{\hangindent}{18pt}

\noindent {\rm (a)} For $z\in\bset$ the bound 
\begin{align}\begin{split}\label{5.0}\big|\big(\tau^z_i(x)\cdot\nabla_x+\tau^z_i(y)\cdot\nabla_y\big)D_x^\alpha D_y^\beta &G_{D,N}(x,y)\big|\\&\le C(\set,k)||x-y|^{2-n-(|\alpha|+|\beta|)}
\end{split}\end{align}
holds for $x,y\in\overline{\sw^z_{R_3}},\,|\alpha|,\,|\beta|\le k_1-1$ and $i=1,\ldots,n-1$. 
\smallskip

\setlength{\hangindent}{18pt} 
\noindent {\rm (b)} Given $z\in\bset,\,R\in (0,R_3]$ and a vector field $\sigma\in C^{0,s}({\sw^z_R})$ for some $s\in (0,1]$ where the extension of $\sigma$ to $\overline{\sw^z_R}$ satisfies $\sigma\perp\nu$ on $\bset\cap\overline{\sw^z_R}$, the bound
\begin{align}\begin{split}\label{5.1}\big|\big(\sigma(x)\cdot\nabla_x+\sigma(y)&\cdot\nabla_y\big)D_x^\alpha D_y^\beta G_{D,N}(x,y)\big|\\
&\le C(\set,k)|\sigma|_{C^{0,s}(\sw^z_R)}|x-y|^{2-n-(|\alpha|+|\beta|)+(s-1)}
\end{split}\end{align}
holds for $x,y\in\overline{\sw^z_R}$ and $|\alpha|,\,|\beta|\le k_1-1$.
\smallskip

	\setlength{\hangindent}{18pt}
\noindent {\rm (c)} For $z\in\bset$ the bound
\begin{equation}\label{5.3}
\big|q\cdot(\nabla_x+\nabla_y\big)D_x^\alpha D_y^\beta G_{D,N}(x,y)\big|\le C(\set,k)|q||x-y|^{2-n-(|\alpha|+|\beta|)}
\end{equation}
holds for $x,y\in\overline{\sw^z_R}$, $|\alpha|,\,|\beta|\le k_1-1$ and vectors $q$ satisfying either $q\perp\nu(x)$ or $q\perp\nu(y)$.
\end{theorem}

\smallskip

\noindent{\bf Proof.}  We can choose $R_3(\set,k)\le R_0$ sufficiently small so that if $z\in\bset,\,R\le R_3$ and $y\in\scy^{z}_R$, then $\db(y)\le \kappa_0/2$ and $\scy^{y'}_R\subset\scy^{z_0}_{R_0}$. This insures that the vector fields $\tau^{z}_i$ and $N^{z}$ are defined in $\scy^{y'}_R$. The constant $R_3$ may be reduced further as the proof proceeds.
\smallskip

\noindent{\bf Step 1: Proof of (a) for $|\alpha|=|\beta|=0$.}

Fix $z_0\in\bset$ and $x_0\in \overline{\sw^{z_0}_{R_3}}$ as in the statement. We will prove \eqref{5.0} for $i=1$ assuming without loss of generality that $x_0\in\set$ and $\db(x_0)\le \kappa_0/2$.  As in the proof of Theorem 4.1 we define the following sets for $R\in(0,R_3]$:
\begin{equation*}
\sm_D(R)=\sup_{y}\{\big|\big(\tau^{z_0}_1(x_0)\cdot\nabla_x+\tau^{z_0}_1(y)\cdot\nabla_y\big)h_D(x_0,y)\big||x_0-y|^{n-2}\}
\end{equation*}
where the sup is over  $y\in B_R(x_0)\cap\overline\set$ with $\db(y)\le \db(x_0)$; and 
$$\sm_N(R)=\sup_{y}\{\big|\big(\tau^{z_0}_1(x_0)\cdot\nabla_x+\tau^{z_0}_1(y)\cdot\nabla_y\big)G_N(x_0,y)\big||x_0-y|^{n-2}\}
$$
where the sup is over  $y\in B_R(x_0)\cap\overline\set$.
We will prove that if $R$ is sufficiently small depending only on $\set$ and $k$ then
\begin{equation}\label{5.4}\sm_{D,N}(R)\le C(\set,k,R^{-1})+C(\set,k)R\sm_{D,N}(R).
\end{equation}
It will then follow that if $R$ is reduced further, again depending only on $\set$ and $k$, then $\sm_{D,N}(R)\le C(\set,k,R^{-1})$ and therefore that \eqref{5.0} holds for points $y_0$ in $B_R(x_0)\cap\overline\set$ (with $\db(y_0)\le \db(x_0)$ in the Dirichlet case). Note that since $(\nabla_x+\nabla_y)\Gamma=0$ the bound for $h_D$ implies the required bound for $G_D$ and in either case we can apply Lemma 4.3(c) for points $y_0\in \overline{\sw^{z_0}_{R_3}}-B_R(x_0)$. The restriction that $\db(y_0)\le \db(x_0)$ for the Dirichlet case is then removed by symmetry. This will prove that \eqref{5.0} holds with $i=1$ and similarly for all $i$.

To prove \eqref{5.4} we fix a point $y_0$ in the defining set of $\sm_{D,N}$ and may assume the following: $y_0\in\set$, again by \eqref{2.115}; that $\db(y_0)\le 
\kappa_0<R_0$ since $\db(x_0)\le\kappa_0/2$ and $R$ is small; that in an appropriate coordinate system 
$y_0'=0$ and $y_0=(0,-\db(y_0))$; and that $\psi^{z_0}$ is transformed in this coordinate system to a $C^k$ function $\psi$ satisfying $\psi(0)=0$ and $\nabla\psi(0)=0$, and the $\tau^{z_0}_i$ are transformed to vector fields $\tau_i=(e_i+\psi_{y_i}e_n)$ and $N^z$ to $(-\psi_{z_1}(\tilde w),\ldots, -\psi_{z_{n-1}}(\tilde w),1)$. In particular, $\tau_i(y_0)=e_i$. Abusing notation slightly we write
$$\sw^{y_0'}_{R_0}=\{(\tilde y, \psi (\tilde y)): |\tilde y|<{R_0}\ {\rm and}\ -{R_0}<y_n<\psi(\tilde y)\}$$ 
and have to prove that
\begin{align}\begin{split}\label{5.5}|\big(\tau_1(x_0)\cdot\nabla_x+\tau_1(y_0)&\cdot\nabla_y\big)(h_D,G_N)(x_0,y_0)||x_0-y_0|^{n-2}\\ &\le C(\set,k,R^{-1})+C(\set,k)R\sm_{D,N}(R).
\end{split}\end{align}
Then taking the sup over $y_0$ in the sets defining $\sm_{D,N}$ we conclude that \eqref{5.4} holds, thus completing Step 1.

Throughout the argument terms bounded by $C(\set,k,R^{-1})+\break C(\set,k)|x_0-y_0|^{2-n}$ will be denoted simply by $O(\se)$. We begin with the Dirichlet case. From Lemma~4.5(a) we have that for $x$ and $y$ near $x_0$ and $y_0$,
\begin{align}\begin{split}\label{5.6}
h_D(x,y)=-\int_\bset& (\dnz \tilde \Gamma(y,z))\Gamma(x,z)dS_z\\
	&-\int_\bset\tilde \Gamma(y,z)\dnz h_D(x,z)dS_z
\end{split} \end{align}
where $\tilde\Gamma(y,z)=\Gamma(y,z)-\Gamma(\ov y,z)$ and $\bar y=(\tilde y,-y_n)$ if $y=(\tilde y,y_n)$ (recall Lemma~4.4(a)). Define
\begin{equation}\label{5.65}
H(y)=\big(\tau_1(x_0)\cdot\nabla_x+\tau_1(y)\cdot\nabla_y\big)h_D(x_0,y).
\end{equation} 
Computing from \eqref{5.6} and noting that $\tau_1(y_0')=\tau_1(0)=e_1$ we find that
$$-H(x_0,y_0)=I+II+III+IV$$
where
\begin{align}\begin{split}\label{5.7}
&I=\int_\bset\big(\dnz\tilde\Gamma(y_0,z)\big)\big(\tau_1(x_0)\cdot\nabla_x\Gamma(x_0,z)\big)dS_z\\
&II=\int_\bset\tilde\Gamma_{y_1,z_j}(y_0,z)\,\Gamma(x_0,z)\nu^j(z)dS_z\\
&III=\int_\bset\tilde\Gamma(y_0,z)\,\dnz\big(\tau_1(x_0)\cdot\nabla_xh_D(x_0,z)\big)dS_z\\
&IV=\int_\bset\tilde\Gamma_{y_1}(y_0,z)\,\dnz h_D(x_0,z)dS_z.
\end{split}\end{align}
We will integrate by parts in $II$ then combine with $I$ to show that $I+II=O(\se)$. First, if $|z-x_0|\ge 2R$ then $|z-y_0|\ge R$ so that the integral in $II$ over the set where $|z-x_0|\ge 2R$ is $O(\se)$. We parameterize the remaining integral by writing $z=(\tilde z,\psi(\tilde z))$ and $\nu(z)dS_z=N(z)d\tilde z $ and noting that $\tilde\Gamma_{y_1,z_j}=-\tilde\Gamma_{z_1,z_j}$. Thus
\begin{align*}\begin{split}
II= O(\se)+\int_{V_{2R}}\Big[&-\frac{\partial}{dz_1}\tilde\Gamma_{z_j}(y_0,(\tilde z,\psi(\tilde z))
	\\&+\tilde\Gamma_{z_j,z_n}(y_0,(\tilde z,\psi(\tilde z))\psi_{z_1}(\tilde z)\Big]\Gamma(x_0,(\tilde z,\psi(\tilde z)) N^j(\tilde z)d\tilde z
		\end{split}\end{align*}
where $V_{2R}=\{\tilde z  :|(\tilde z,\psi(\tilde z))-x_0|<2R\}$. Integrating by parts with respect to $z_1$ and checking that the boundary term is $O(\se)$ we obtain
\begin{align}\begin{split}\label{5.8}
II= O(\se)+\int_{V_{2R}}\Big(\big[&\nabla_z\tilde\Gamma(y_0,z)\cdot N(z)\big]\big[\tau_1(z)\cdot\nabla_z\Gamma(x_0,z)\big]\\
	&\quad\qquad\qquad+\tilde\Gamma_{z_j}(z)\Gamma(x_0,z)N^j_{z_1}(z)\\
		&+\tilde\Gamma_{z_j,z_n}(y_0,z)\Gamma(x_0,z)\psi_{z_1}(\tilde z)N^j(z)\Big)d\tilde z
\end{split}\end{align}
where $z=(\tilde z,\psi(\tilde z))$. We replace $\tau_1(z)$ by $\tau_1(x_0)$ on the right, incurring a change in the integrand which by Lemma 4.5(b) is bounded by $C(\set,k)\db(y_0)|z-y_0|^{-n}|z-x_0|^{2-n}$, and therefore an $O(\se)$ change in the integral. (To see this let $\delta=2|x_0-y_0|/3$ and consider the integrals over $\bset\cap B_\delta(x_0),\,\bset\cap B_\delta(y_0)$ and the complement of their union. The argument for $\bset\cap B_\delta(y_0)$ requires explicit use of the term $\db(y_0)$ appearing in Lemma~4.5(b).) In the second term on the right we have that $j\not=n$ so that the integrand is $O(\se)$, again by Lemma 4.5(b), with the same result. The same is true for the third term in the integral on the right for $j\not= n$ because $|\nabla_z\psi(\tilde z)|\le C|\tilde z|\le C|z-y_0|$ by Lemma~4.4(b). For $j=n$ we write $\tilde\Gamma_{z_n,z_n}=-\tilde\Gamma_{z_1,z_1}+\ldots$ and integrate by parts with respect to $z_1,\ldots$. The resulting integrals are then bounded as for the previous terms. Thus 
\begin{align}\label{5.9}
II=O(\se)	+\int_{\bset\cap B_{2R}(x_0)}\dn\tilde\Gamma(y_0,z)\big(\tau_1(x_0)\cdot\nabla_z\Gamma(x_0,z)\big)dS_z.
\end{align}
 We can also replace the domain of integration in $I$ by $\bset\cap B_{2R}(x_0)$, just as for $II$. The resulting integral then cancels the integral in \eqref{5.9} because $(\nabla_x+\nabla_z)\Gamma=0$. We conclude that
\begin{equation}\label{5.10}
|I+II|\le C(\set,k,R^{-1})+C(\set,k)|x_0-y_0|^{2-n}=O(\se).
\end{equation}

Next we estimate term $IV$ and combine the result with $III$. First, just as for $II$ above, the integral in $IV$ over the set where $|z-x_0|\ge 2R$ is $O(\se)$ by Lemma~4.3(c) (it is at this step that we require $k_1\ge 1$).  The integral over the remaining set can again be parameterized by writing $z=(\tilde z,\psi(\tilde z))$ and $\nu(z)dS_z=N(z)d\tilde z $. In the following computation we integrate by parts with respect to $z_1$, checking that the boundary term is $O(\se)$, omitting the arguments $(y_0,z)$, $(x_0,z)$ and $\tilde z$ on the $\tilde\Gamma,\,\Gamma$ and $N$ terms respectively and writing $h$ in place of $h_D$:

\begin{align*}
IV=O(\se) +\int_{V_{2R}}\Big[-\frac{\partial}{\partial z_1}\tilde\Gamma(y_0,z)+\tilde\Gamma_{z_n}(y_0,z)\psi_{z_1}(\tilde z)\Big]h_{z_j}(x_0,z)N^j(\tilde z)d\tilde z\\
=O(\se) +\int_{V_{2R}}
\Big[\tilde\Gamma(h_{z_1,z_j}+\Gamma h_{z_n,z_j}\psi_{z_1})N^j+\tilde\Gamma h_{z_j}N^j_{z_1}+\tilde\Gamma_{z_n}h_{z_j}\psi_{z_1}N^j\Big]d\tilde z.\end{align*}
The term in parentheses inside the brackets is
$$\big(h_{z_1}+h_{z_n}\psi_{z_1}\big)_{z_j}-h_{z_n}\psi_{z_1,z_j}=\big[\tau_1(z)\cdot\nabla_z h(x_0,z)\big]_{z_j}-h_{z_n}\psi_{z_1,z_j}$$
so that
\begin{align*}
IV=O(\se) + \int_{\bset\cap B_{2R}(x_0)}\big[\tilde\Gamma(y_0,z)\dnz\big(\tau_1(z)\cdot\nabla_z h_D(x_0,z)
\big)dS_z\\
+\int_{V_{2R}}\big[-\tilde\Gamma h_{z_n}\psi_{z_1,z_j}N^j+\tilde\Gamma h_{z_j}N^j_{z_1}+\tilde\Gamma_{z_n}h_{z_j}\psi_{z_1}N^j\big]d\tilde z.\end{align*}
The second integral on the right here is $O(\se)$ just as for term $II$ because by \eqref{2.11} and Lemma~4.5(b) $h_D$ satisfies the same pointwise bounds \eqref{1.1} satisfied by $\Gamma$ (the assumption that $\db(y_0)\le \db(x_0)$ is applied here). Combining with the integral in $III$ and applying \eqref{5.10} we therefore have that 
$$H(y_0)=O(\se)- \int_{\bset\cap B_{R/2}(x_0)} \tilde\Gamma(y_0,z)\dnz H(z)dS_z,$$
the change in the domain of integration incurring an error bounded by $O(\se)$. Next we
define a harmonic function $H'(y)$ to be the same as $H(y)$ in \eqref{5.65} but with $\tau_1(y)$ replaced by $\tau_1(x_0)$. Then since  $H'(y_0)-H(y_0)=O(\se)$ by Theorem 4.1, $H(z)$ can be replaced by $H'(z)$ in the integral on the right above at the expense of another $O(\se)$ term. Thus
\begin{equation}\label{5.165}
H'(y_0)=O(\se)- \int_{\bset\cap B_{R/2}(x_0)} \tilde\Gamma(y_0,z)\dnz H'(z)dS_z.
\end{equation} 
Now if we define 
\begin{equation}\label{5.166}\sm'_D(R)=\sup_{y}\{|x_0-y|^{n-2}|H'(y)|\}
\end{equation}
where the sup is over the same set as in the definition of $M_D(R)$ in \eqref{2.175} then \eqref{5.166} is the same as \eqref{2.175} but with $H'$ in place of $h_D$; and \eqref{5.165} implies the bound in \eqref{2.19}, again with $H'$ in place of $h_D$. Since $H'$ is harmonic, we can therefore repeat the argument proceeding from \eqref{2.19} to find that
$$|x_0-y_0|^{n-2}| H'(y_0)|\le C(\set,k,R^{-1})+C(\set,k)R\sm'_D(R)$$
just as in \eqref{2.195}. The same bound then holds for  $H(y_0)$, and this proves \eqref{5.5} for the Dirichlet case.

Next we prove \eqref{5.5} for the Neumann case. From Lemma~4.5(a) we have that

\begin{align}\begin{split}\label{5.11}
\tilde G_N(x,y)=\tilde\Gamma(x,y)-&\ov{\Gamma(\cdot,y)} + \ov {\Gamma(\cdot,\ov y)}\\
	&+\int_\bset\big(\dnz \tilde\Gamma(y,z)\big)\tilde G_N(x,z)dS_z
\end{split}\end{align}
for $x$ and $y$ near $x_0$ and $y_0$, where $\tilde G_N=\Gamma + h_N$ and now $\tilde\Gamma(y,z)=\Gamma(y,z)+\Gamma(\ov y,z)$. We define a new function $H$ by
\begin{equation}\label{5.115}
 \quad\qquad H(y)=\big(\tau_1(x_0)\cdot\nabla_x+\tau_1(y)\cdot\nabla_y\big)\tilde G_N(x_0,y)
\end{equation}
and compute from \eqref{5.11} that
$$H(y_0)=I+II+III$$
where $I$ is the result of applying the operator $\tau_1(x_0)\cdot\nabla_x+\tau_1(y_0)\cdot\nabla_y$  to the three terms preceding the integral on the right in \eqref{5.11} and 
$$\qquad II=\int_\bset \big(\dnz\tilde\Gamma(y_0,z)\big)\big(\tau(x_0')\cdot\nabla_x\tilde G_N(x_0,z)\big)dS_z$$
and 
$$III=\int_\bset\tilde\Gamma_{z_j,y_1}(y_0,z)\tilde G_N(x_0,z)\,\nu^j(z)\,dS_z.$$
It is easy to see that $I= O(\se)$ and we can restrict the domain of integration in both $II$ and $III$ to $\bset\cap B_{2R}(x_0)$, the integrals in the complement being $O(\se)$. We parameterize the remaining part of $III$ by writing $z=(\tilde z,\psi(\tilde z))$, $\nu\, dS=Nd\tilde z$ and $\tilde\Gamma_{y_1,z_j}=-\tilde\Gamma_{z_1,z_j}$. Integrating by parts with respect to $z_1$ and checking that the boundary term is $O(\se)$ we obtain
\begin{align}\begin{split}\label{5.12}
III=&O(\se) + \int_{V_{2R}}\Big[-\frac{\partial}{\partial z_1}\tilde\Gamma_{z_j} (y_0,z)+\tilde\Gamma_{z_j,z_n}(y_0,z)\psi_{z_1}(\tilde z)\Big]\tilde G(x_0,z)N^j(\tilde z)d\tilde z\\
=&O(\se) +\int_{\bset\cap B_{2R}(x_0)}\big(\dnz \tilde\Gamma(y_0,z)\big)(\tau_1(z)\cdot\nabla_z) \tilde G(x_0,z)\,dS_z + IV
\end{split}\end{align}
where
$$IV=\int_{V_{2R}}\big[\tilde\Gamma_{z_j}(y_0,z) N^j_{z_1}(\tilde z)+\tilde\Gamma_{z_j,z_n}(y_0,z)\psi_{z_1}(\tilde z)N^j(\tilde z)\big]\tilde G(x_0,z)\, d\tilde z.$$
Combining this with $II$ we then have that
\begin{equation}\label{5.13}
H(y_0)=O(\se) +\int_{\bset\cap B_{2R}(x_0)} \big(\dnz\tilde\Gamma(y_0,z))\,H(z)dS_z + IV.
\end{equation}
Term $IV$ is computed in the following steps: we substitute $N^j=-\psi_{z_j}$ for $j<n$ and $N^n=1$, rearrange terms, integrate by parts, noting that the boundary integral over $\partial V_{2R}$ is $O(\se)$, and identify divergence free vector fields $K^j=-\tilde G_{z_n}e_j+\tilde G_{z_j}e_n$ for $j<n$. The result is that

\begin{align*}
IV&=-\int_{V_{2R}}\sum_{j<n}\frac{\partial}{\partial z_j}\big[\psi_{z_1}(\tilde z)\tilde\Gamma_{z_j}(y_0,z)\big]\tilde G(x_0,z)d\tilde z\\
&=O(\se) + \sum_{j<n}\int_{\bset\cap B_{2R}(x_0)}\psi_{z_1}(\tilde z)\tilde\Gamma_{z_j}(y_0,z)K^j(x_0,z)\cdot \nu(z)dS_z.
\end{align*}
To bound the integral here for fixed $j<n$ we  replace the domain of integration by $\partial \set^\ep_{2R}$ where
$$\set^\ep_{2R}=(\set\cap B_{2R}(x_0))-B_\ep(x_0)-B_\ep(y_0)$$ 
for small $\ep>0$ then take $\lim_{\ep\to 0}$. To justify this we let $\mathcal{I}(z)$ be the integrand and note first that $\mathcal{I}$ is $O(\se)$ on $ \partial B_{2R}(x_0)$ hence so is its integral over $ \partial B_{2R}(x_0)$; and we claim further that
$$\lim_{\ep\to 0}\int_{\partial B_\ep(x_0)} \mathcal{I}(z)dS_z= \lim_{\ep\to 0}\int_{\partial B_\ep(y_0)} \mathcal{I}(z)dS_z=0.$$
To prove these statements we can first replace $\tilde G_N$ by $\Gamma$ in the definition of $K^j$, hence of $\mathcal{I}$, because $h_N(x_0,\cdot)$ is bounded independently of $\ep$ by Theorem 3.1(a). The remaining terms in $\mathcal{I}$ are then identically zero on $\partial B_\ep(x_0)$ by direct computation; and for $z\in \partial B_\ep(y_0)$, $|\mathcal{I}(z)|\le C|z-y_0|^{2-n}|x_0-y_0|^{1-n}$, whose integral goes to zero as $\ep\to 0$. We therefore have that
\begin{align*}
IV=&O(\se)  +\lim_{\ep\to 0}\sum_{j<n} \int_{\partial \set^\ep_{2R}}\psi_{z_1}(\tilde z)\tilde\Gamma_{z_j}(y_0,z)K^j(x_0,z)\cdot \nu(z)dS_z\\
&=O(\se)  + \lim_{\ep\to 0}\int_{\set^\ep_{2R}} \nabla_z\big(\psi_{z_1}(\tilde z)\tilde\Gamma_{z_j}(y_0,z)\big)\cdot K^j(x_0,z)\,dz
\end{align*}
because $K^j(x_0,\cdot)$ is divergence free. The integrand on the right here is bounded by $C(|z-y_0||z-x_0|)^{1-n}$ so that the limit is $O(\se)$, and therefore  $IV$ is $O(\se)$. Returning to \eqref{5.13}  we thus obtain 
$$H(y_0)=O(\se)+\int_{\bset\cap B_{2R}(x_0)} \big(\dnz\tilde\Gamma(y_0,z))H(z)dS_z .$$
We can replace the domain of integration here by $\bset\cap B_{R/2}(x_0)$ at the expense of another $O(\se)$ term and then apply Lemma 4.5(b) to bound the integrand by $C(\set,k)\sm_N(R)(|z-y_0||z-x_0|)^{2-n}$. The integral is then bounded by $C(\set,k)\sm_N(R)R|x_0-y_0|^{2-n}$ just as in the discussion preceding \eqref{2.161}. We then multiply by $|x_0-y_0|^{n-2}$ and take the sup over $y_0$ in the set defining $\sm_N(R)$ to obtain finally that
$$\sm_N(R)\le C(\set,k,R^{-1})+C(\set,k)R\sm_NR).$$
This proves \eqref{5.4} and completes the proof of Step 1 for the Neumann case. \medskip

\noindent{\bf Step 2: Proof of (a) for $|\alpha|,\,|\beta|\le k_1-1$.} Fix $z_0\in\bset$ and $\nobreak{x_0,\,y_0\in \sw^{z_0}_{R_3}}$ and let $\delta=|x_0-y_0|$. We may assume that $\delta\le R_1$, otherwise Lemma 4.3(c) applies to show that the left side of \eqref{5.0} is bounded by $C(\set,k)$ since $|\alpha|+1,\,|\beta|+1\le k_1$. Recalling Lemma 4.3(b), we therefore have that $s(x_0,y_0)=\delta$ and so in the nontrivial case of that lemma that $\db(y)<\delta/12,$
$${\rm dist}(x_0,\sw^{y'_0}_{\delta/3})\ge C^{-1}\delta\ {\rm and}\ y_0\in  \sw^{y'_0}_{\delta/6}.$$ 
It then follows from \eqref{2.12} in Theorem 4.1 that for $|\alpha|\le k_1$ and $|\gamma|\le k$
\begin{equation}\label{5.19}
|D_x^\alpha D_y^\gamma G_{D,N}(x_0,\cdot)|_{L^2(\sw^{y'_0}_{\delta/3})}\le C\delta^{2-n/2-(|\alpha|+|\gamma|)}.
\end{equation}
We are going to apply Theorem 7.2(b) to the function 
\begin{equation}\label{5.195}
K_\alpha(y)\equiv \big(\tau_i(x_0)\cdot\nabla_x+\tau_i(y)\cdot\nabla_y\big)D_x^\alpha G_{D,N}(x_0,y)
\end{equation}
for fixed $\alpha$ with $|\alpha|\le k_1-1$. To facilitate the computations we denote by $\se^{p,q}_i(y)$ generic linear combinations of derivatives $D_x^\zeta D_y^\eta G_{D,N}(x_0,y)$ for $|\zeta|\le p$ and $|\eta|\le q$ and with coefficients which are $C^{k-i}$ functions with derivatives bounded by $C(\set,k)$ (the latter will be determined by $i$-th order derivatives of the $\tau_j$). We check that $\Delta K_\alpha=\se^{|\alpha|,2}_2$ on $\sw^{y'_0}_{\delta/3}$ so that in the Dirichlet case
$$\int_{\sw^{y_0'}_{\delta/3}}\nabla K_\alpha\cdot\nabla v\, dy=-\int_{\sw^{y_0'}_{\delta/3}} \Delta_y K_\alpha\,v\, dy = \int_{\sw^{y_0'}_{\delta/3}} \se^{|\alpha|,2}_2v\, dy$$
for $v\in H^1_0(\sw^{y_0'}_{\delta/3})$. Note also that $K_\alpha =0$ on the top of $\sw^{y_0'}_{\delta/3}$, as required in case (D) of Theorem 7.2. In the Neumann case we take $v\in H^1(\sw^{y_0'}_{\delta/3})$ vanishing on the sides and bottom of $\sw^{y_0'}_{\delta/3}$ and compute that 
$$\int_{\sw^{y_0'}_{\delta/3}}\nabla K_\alpha\cdot\nabla v=\int_{\partial\sw^{y_0'}_{\delta/3}}(\dn K_\alpha)\,v-\int_{\sw^{y_0'}_{\delta/3}} \se^{|\alpha|,2}_2v. $$
Now $\dn K_\alpha =\se^{|\alpha|,1}_1$ on $\bset$ because $\dn G_N(x_0,\cdot)=0$. Applying the trace theorem to the boundary integral, we therefore conclude that
$$\int_{\sw^{y_0'}_{\delta/3}}\nabla K_\alpha\cdot\nabla v=\int_{\sw^{y_0'}_{\delta/3}} \big(\se^{|\alpha|,2}_2v +\se^{|\alpha|,1}_1\,\nabla v\big)$$
in either case\ (D) or (N) of Theorem 7.2. We now apply \eqref{a.32} with $u=K_\alpha,\,\varphi=\se^{|\alpha|,2}_2$ and $\Phi =\se^{|\alpha|,1}_1$ and note that
for $i\le j-1\le k-1$,
$$|\se^{|\alpha|,2}_2|_{H^{i-1}},\  |\se^{|\alpha|,1}_1|_{H^i}= |\se^{|\alpha|,i+1}_{i+1}|_{L^2} \le C\delta^{2-n/2-(|\alpha| +i+1)} $$
(norms are over $\sw^{y_0'}_{\delta/3}$) by \eqref{5.19}. The result is that for $|\gamma|\le k$,
\begin{equation}\label{5.20}
|D_y^\gamma K_\alpha|_{L^2(\sw^{y_0'}_{\delta/6})} \le C(\set,k)\big[\delta^{-|\gamma|}|K_\alpha|_{L^2(\sw^{y_0'}_{\delta/3})} + \delta^{2-n/2-(|\alpha|+|\gamma|)}\big].
\end{equation}

Now for the case that $\alpha=0$ we can apply the pointwise bound proved in Step~1 to compute that $|K_0|_{L^2(\sw^{y_0'}_{\delta/3})}\le C\delta^{2-n/2}$ and therefore that 
$$|D_y^\gamma K_0|_{L^2(\sw^{y_0'}_{\delta/6})} \le C(\set,k)\delta^{2-n/2-|\gamma|}$$
for $|\gamma|\le k$. We then apply the scaling result Theorem 7.4(a) just as in \eqref{2.20}--\eqref{2.21} in the proof of Theorem 4.1 to find that 
$$\sup_{\sw^{y_0'}_{\delta/6}}|D_y^\beta K_0| \le C(\set,k)\delta^{2-n/2-|\beta|}$$
for $|\beta|\le k_1$. Thus in particular
$$\Big|D^\beta_y\big[\big(\tau_1(x)\cdot\nabla_x+\tau_1(y)\cdot\nabla_y\big)G_{D,N}(x,y)\big]\Big|\le C(\set,k)|x-y|^{2-n-|\beta|}$$
at $(x_0,y_0)$ and therefore for all $x,\,y\in \sw^{z_0}_{R_3}$ and $|\beta|\le k_1-1$ (which was assumed earlier in Step~2). The operator $D_y^\beta$ can then be passed under the directional derivatives at the expense of an error which by \eqref{2.11} can be absorbed into the term on the right. 
Applying symmetry and replacing $\beta$ by $\alpha$ we conclude that
\begin{equation}\label{5.205}|(\tau_1(x)\cdot\nabla_x+\tau_1(y)\cdot\nabla_y\big)D_x^\alpha G_{D,N}(x,y)|\le C(\set,k)|x-y|^{2-n-|\alpha|}
\end{equation}
for $|\alpha|\le k_1-1$ and $x,\,y\in \sw^{z_0}_{R_3}$. 

We now return to the points $x_0$ and $y_0$ fixed at the beginning of Step 2, again with $\delta=|x_0-y_0|$ and $K_\alpha$ as in  \eqref{5.195}. The left side of \eqref{5.205} with $x=x_0$ is precisely  $|K_\alpha(y)|$ so that $|K_\alpha|_{L^2(\sw^{y_0'}_{\delta/3})}\le C\delta^{2-n/2-|\alpha|}$. Substituting this on the right in \eqref{5.20} we then find that $|D_y^\gamma K_\alpha|_{L^2(\sw^{y_0'}_{\delta/6})} \le C(\set,k)\delta^{2-n/2-(|\alpha|+|\gamma|)}$ for $|\gamma|\le k$. Then just as for the $K_0$ case we can apply Theorem 7.4(b) with the result that \eqref{5.0} holds for $|\alpha|,\,|\beta|\le k_1-1$ and points $x,\,y\in \sw^z_{R_3}$, for all $z\in\bset$. This completes the proof of part (a) of the theorem.

\bigskip
\noindent{\bf Step 3: Proofs of (b) and (c).} We can replace the term $\tau^z_i(y)$ on the left side of \eqref{5.0} by $\tau^z_i(x)$ at the expense of a term bounded by $\nobreak{C|x-y||\nabla_y D_x^\alpha D_y^\beta G_{D,N}(x,y)|}$. Since $\nobreak{|\beta|+1\le k_1}$, \eqref{2.11} applies to show that the latter is bounded by $C(\set,k)|x-y|^{2-n-(|\alpha|+|\beta|)}$, which can be absorbed into the right side of \eqref{5.0}. The same is true if $\tau^z_i(x)$ is replaced by $\tau^z_i(x')$, again by \eqref{2.11} and Corollary 4.2. We thus have that for $z\in \bset$ and $x,y\in \sw^z_{R_3}$,
\begin{align}\begin{split}\label{5.21}\big|\tau^z_i(x')\cdot\big(\nabla_x+\nabla_y\big)D_x^\alpha D_y^\beta &G_{D,N}(x,y)\big|\\&\le C(\set,k)|x-y|^{2-n-(|\alpha|+|\beta|)}.
\end{split}\end{align}
Note that this reduction required only the Lipschitz regularity of $\tau^z_i$. Now if $q\perp \nu(x')$ then $q$ is a linear combination of $\tau^z_1(x'),\ldots,\tau_{n-1}^z(x')$ and one checks easily that the coefficients in this linear combination are bounded by $C(\set,k)|q|$. Thus  \eqref{5.3} holds, and this proves (c). Similarly, if $\sigma$ is a vector field as in (b), then \eqref{5.21} holds with $\sigma(x')$ in place of $q$. The steps leading from \eqref{5.0} to \eqref{5.21} with $\sigma$ in place of $\tau^z_i$ can then be reversed (with the H\"older-$s$ regularity of $\sigma$ in place of the Lipschitz regularity of the $\tau_i$) to prove \eqref{5.1}. 

\rightline{\qedsymbol}

\bigskip

		\section{An Application to Fluid Mechanics}
\bigskip
In this section we describe the application of the pointwise bounds in Theorems 4.1 and 5.1 to the construction of solutions of the Navier-Stokes equations of barotropic, compressible flow on fairly general subsets of $\R^3$. This construction extends earlier results in \cite{halfspace} for flows in $\R^3$ and in half-spaces of $\R^3$, now enabled in particular by the cancellation property in Theorem 5.1(b) for the Green's function for the Neumann-Laplace operator. Complete details will be reported elsewhere. Here we will content ourselves with extracting a single {\it a priori} bound required in the proof and describing informally how the results of the present paper are applied in an essential way.

The Navier-Stokes equations describe a compressible fluid occupying a bounded open set $\set$ in $\R^3$ in terms of its density $\rho$ and the particle velocity 
$u=(u^1,u^2,u^3)$ which are the unknown functions of $x\in\set$ and time $t\ge 0$. The evolution of $\rho$ is governed by the principle of conservation of mass, expressed by the familiar equation
\begin{equation}\label{6.1}
\rho_t+{\rm div}(\rho u)=0;
\end{equation}
and the evolution of the velocity is governed by Newton's law, as follows. If $x(t)$ is the position at time $t$ of a fixed fluid particle, then its velocity is $u(x(t),t)$ and therefore the $j$-th component of its acceleration is $\dot u^j\equiv u^j_t+(\nabla_xu^j)\cdot u$. Thus the mass (per volume) of a fluid particle times its acceleration is 
$$\rho\dot u^j=\rho (u^j_t +\nabla_xu^j\cdot u)=(\rho u^j)_t+(\rho u^ju^k)_{x_k}$$
(sum over $k$), and this equals the sum of the forces acting on the particle. We need not describe these forces in detail but instead simply note that for each fixed $t$ the sum of the forces is a vector field on $\set$ which can be expressed as the gradient of a scalar potential $F$ plus a divergence-free vector field $H$. Thus
\begin{equation}\label{6.2}
\rho\dot u=\nabla_x F + H.
\end{equation}
In the Navier-Stokes system $F$ and $H$ are determined by $\rho,\,u$ and $\nabla u$ and can therefore  be regarded as functions of $x$ and $t$ by composition. 
A boundary condition is imposed on the velocity and this includes the condition $u(x,t)\cdot\nu(x)=0$ for $x\in\bset$. Initial data $(\rho_0,u_0)$ is given and the system \eqref{6.1}-\eqref{6.2} together with the boundary and initial conditions are to be solved for $t\ge 0$ under fairly mild restrictions on the system parameters and under the assumptions  that $\rho_0$ is measurable and takes values in a compact subset of $(0,\infty)$, and that $u_0\in L^6(\set)$ and the initial energy, kinetic plus potential, satisfies a certain size restriction. 

The proof consists in a set of {\it a priori} bounds for smooth approximate solutions, bounds which are fairly technical and which are coupled in complicated ways. A crucial step in this construction is to bound the quantity
${\mathcal F}\equiv\displaystyle\int_0^{\bar t} F(x(t),t)\,dt$, where $x(t)$ is a particle trajectory near $\bset$, in terms of a bound $C_1$ for the quantities
\begin{equation}\label{6.3}\sup_t|\rho(\cdot,t)|_{L^\infty},\,\sup_t|u(\cdot,t)|_{L^6}\,\ {\rm and}\ \int_0^{\bar t}|u(\cdot,t)|_{C^{0,s}}|u(\cdot,t)|_{L^q}dt \end{equation}
(norms over $\set$) for particular $s\in(0,1)$ and $q>6$ satisfying $\nobreak{q'(3-s)<3}$ determined elsewhere in the proof. We  will describe how the required bound for ${\mathcal F}$ is obtained for an approximate solution $(\rho,u)$ defined up to time $\bar t$ where $\rho,\,u$ and $\bset$ are as smooth as required for the arguments that follow. 

First we compute from \eqref{6.2} that at each $t$,
\begin{align}\begin{split}\label{6.4}
\begin{cases}\Delta_x F = {\rm div}_x(\rho\dot u)\ {\rm in}\ \set\\
\dn F=(\rho\dot u-H)\cdot\nu \ {\rm on}\ \bset.
\end{cases}
\end{split}\end{align}
Observe that the term $\rho \dot u$ here indicates a derivative along a particle path whereas $\mathcal F$ is an integral over a particle path. Cancellation can therefore be anticipated and will be detected through the representation \eqref{2.5} of $F(x(t),t)$ in terms of the Neumann-Green's function $G$:
\begin{align*}
\mathcal{F}=\int_0^{\bar t}&\Big(\overline{F(\cdot,t)}-\int_\bset  G(x(t),y)\dn F(y,t)dS_y \\
	&+\int_\set G(x(t),y)\big[(\rho u^j)_{y_j,t}+(\rho u^ju^k)_{y_j,y_k}\big]\,dy\Big)dt.
\end{align*}
The term involving $\bar F$ is of lower order and together with other such terms will be denoted by an ellipsis. We integrate by parts in the integral over $\set$ with respect to $t$ and $y_k$, noting that first-order derivatives of $G(x(t),\cdot)$ are integrable by \eqref{2.11}.  Applying the relation $\dot x^k(t)=u^k(x(t),t)$ and rearranging we find that, modulo the ellipsis,
\begin{align}\begin{split}\label{6.5}
\mathcal{F}=&\int_0^{\bar t}\!\!\int_\bset G(x(t),y)\big[-\dn F(y,t)+(\rho u^ju^k)_{y_j}(y,t)\nu^k(y)\big]dS_ydt\\
&+\int_\set G(x(\cdot),y)(\rho u^j)_{y_j}(y,\cdot)dy\,\Big|_0^{\bar t}\\
	&-\int_0^{\bar t}\!\!\int_\set \big[G_{x_k}(x(t),y)u^k(x(t),t)(\rho u^j)_{y_j}(y,t)\\
		&\qquad\qquad\qquad\qquad\qquad+G_{y_k}(x(t),y)(\rho u^ju^k)_{y_j}(y,t)\big]dy\,dt.
\end{split}\end{align}
Applying the relations $u\cdot\nu=u_t\cdot\nu=0$ on $\bset$ we find that the term in brackets in the boundary integral here can be written
$$\big[-\nabla F+(\rho u^j)_{y_j}u+\rho(\nabla u\,u)\big]\cdot \nu=\big(-\nabla F +\rho \dot u\big)\cdot\nu=-H\cdot\nu.
$$
This boundary integral thus depends on details of the divergence-free field $H$ and on the specific boundary conditions. The required bound is technical and not germane to the present discussion and so will be omitted. Concerning the single integral on the right in \eqref{6.5} we integrate by parts and then note that $G_{y_j}(x(\cdot),\cdot)$ is in $L^p(\set)$ for $p<3/2$ by \eqref{2.11} and that $(\rho u)(\cdot,t)\in L^6(\set)$ by our assumptions in \eqref{6.3}. Bounds for the first two integrals on the right in \eqref{6.5} can therefore be included in the ellipsis.

Integration by parts appears to be required in the third integral on the right in \eqref{6.5} because derivatives  of $\rho$ and $u$ are not included in \eqref{6.3}. Second derivatives of $G(x,\cdot)$ are not in general integrable, however, so that if this analysis is to succeed, some cancellation in the integrand  must occur. We write
\begin{align}\begin{split}\label{6.6}\!\!\!\!\!\!\mathcal{F}=\ldots-\lim_{\ep\to 0} \int_0^{\bar t}\!\!\int_{\set-B_{\ep}(x(t))} \Big(&G_{x_k}(x(t),y)\big[(\rho u^j)(y,t)u^k(x(t),t)\big]_{y_j}\\&+G_{y_k}(x(t),y)\big[(\rho u^ju^k)(y,t)\big]_{y_j}\Big)dy\,dt.
\end{split}\end{align}
We now integrate by parts with respect to $y_j$ in both terms and note that the boundary integral over $\bset$ is zero because $u(x,\cdot)\cdot\nu=0$ on $\bset$. Concerning the integral over $\partial B_\ep(x(t))$ we recall that by Theorem~3.1 $G(x,\cdot)=\Gamma(x,\cdot)+h(x,\cdot)$ where $h(x(t),\cdot)\in H^{k+1}(\set)$ . The contribution from $h$ to the integral over $\partial B_\ep(x(t))$ therefore vanishes in the limit as $\ep\to 0$, and we need only consider the contribution from $\Gamma$. Modulo multiplicative constants this contribution is
$$\int_0^{\bar t}\!\!\int_{B_\ep(x(t))}\Big(\frac{x_k-y_k}{|x-y|^3}u^k(x(t),t)+\frac{y_k-x_k}{|x-y|^3}u^k(y,t)\Big)(\rho u)(y,t)\cdot\nu(y)\,dS_y\,dt$$
which is easily seen to go to zero as $\ep\to 0$ (recall that $u$ is an approximate smooth solution here). 
We thus have that
\begin{align}\begin{split}\label{6.7}\mathcal{F}&=\ldots+\lim_{\ep\to 0} \int_0^{\bar t}\!\!\int_{\set-B_{\ep}(x(t))} \Big(G_{x_k,y_j}(\rho u^j)(y,t)u^k(x(t),t)\\
&\qquad\qquad\qquad\qquad\qquad\qquad\qquad+G_{y_k,y_j}(\rho u^ju^k)(y,t)\Big)dy\,dt\\
&=\ldots+\lim_{\ep\to 0} \int_0^{\bar t}\!\!\int_{\set-B_{\ep}(x(t))} \Big(\big[G_{x_k,y_j}u^k(x(t),t)+G_{y_k,y_j}u^k(y,t)\big]\Big)\\
&\quad\qquad\qquad\qquad\qquad\qquad\qquad\qquad\qquad\qquad\times(\rho u^j)(y,t)dy\,dt
\end{split}\end{align}
where each term involving $G$ is evaluated at $(x(t),y)$. Contributions to this integral from points $y$ far from $\bset$ or far from $x(t)$ can be absorbed into the ellipsis by Corollary 4.2 (because $G$ can be replaced by $h_N$ in the integrand) or by Lemma~4.3(c). For $x(t)$ and $y$ close and close to $\bset$, Theorem 5.1(b) applies with $|\alpha|+|\beta|=1$ to show that the term in parentheses here is bounded by 
$C|u(\cdot,t)|_{C^{0,s}}|y-x(t)|^{s-3}$. Taking the $L^{q'}$ norm of this term and the $L^q$ norm of $(\rho u^j)(\cdot,t)$ we can then bound the integral in \eqref{6.7} by
$$ C_1\int_0^{\bar t}|u(\cdot,t)|_{C^{0,s}}|u(\cdot,t)|_{L^q}dt$$
since $q'(3-s)<3$. This competes the required bound for $\mathcal{F}$.
\bigskip

\break

\section{Appendix: Interior-type Estimates Near the Boundary} $\ $  
 
 \bigskip
 
In this section we derive regularity results for weak solutions of second-order elliptic equations near the boundary of the spatial domain. Specifically, we consider a solution $u$ of the elliptic problem \begin{equation}\label{a.1}
\int_{\sw_{R_0}}(A\nabla u)\cdot\nabla v= F\cdot v
\end{equation} 
which is to hold for test functions $v$ in a specified space and for elements $F$ in its dual.
Here $A$ is an $n\times n$ matrix-valued function on a spatial domain $\sw_{R_0}\subset\R^n$ described as follows, consistent with the notations in Definition 2.1:\medskip

\noindent\emph{Let  $\mathcal{C}_{R}$ be the cylinder $\mathcal{C}_{R}\equiv B_{R}\times (-R,R)$ where for $R>0$ $B_{R}$ is the ball of radius $R$ centered at the origin in $\R^{n-1}$ where $n\ge 2$. Fix $R_0>0$ and $k \ge 1$ and assume that there is a mapping $\psi\in C^{k+1}(B_{R_0})$ such that $\psi(0)=0$ and $\nabla\psi(0)=0$, that $\psi(B_R)\subset (-R/2,R/2)$ for $R\le R_0$ and that the derivatives of $\psi$ up to order $k+1$ are bounded by a constant $M_\psi(k)$. Then for $R\in (0,R_0]$ define
 $$\sw_R=\{ (\tilde z,z_n) : \tilde z\in B_R {\rm\ and\ } -R<z_n<\psi(\tilde z)\},$$
its top
$$\partial\sw_R^+=\{ (\tilde z,z_n) : \tilde z\in B_R {\rm\ and\ } z_n=\psi(\tilde z)\}$$
and a $C^k$ vector field  $N(z)=(-\nabla\psi(\tilde z),1)$ at points $z=(\tilde z,z_n)$ of $\mathcal{C}_{R_0}$. Thus $\nu\equiv N/|N|$ is the unit outer normal to $\sw_{R_0}$ at points on $\partial\sw_{R_0}^+$. Concerning $A$ we assume that its derivatives up to order $k$ exist and are bounded on $B_{R_0}$ by a constant $M_A(k)$ and that there is a positive constant $\eta_A$ such that for all $x\in \sw_{R_0}$ and all $\xi\in\R^n$}
\begin{equation}\label{a.2}
(A(x)\,\xi)\cdot\xi \ge\eta_A|\xi|^2.
\end{equation}

\medskip

The major result of this section is Theorem 7.2 below in which we show that for small $\delta,\,\delta'\in (0,R_0)$ with $\delta'<\delta$, higher-order regularity norms of a solution $u$ of \eqref{a.1} on the smaller set $\sw_{\delta'}$ are bounded by lower-order regularity norms on the larger set $\sw_{\delta}$ modulo constants and inverse powers of $\delta-\delta'$. Sobolev norms are considered in Theorem 7.2(b) and Lebesgue norms in Theorem 7.2(c). The proofs apply standard cutoff functions which modulate the two sets in question but which must also respect the natural boundary condition associated to \eqref{a.1} on $\partial\sw_\delta^+$:

\begin{lemma} Assume that the hypotheses and notations above are in force. Then given $\theta\in (0,1)$ there is a number $\delta_0\in (0,R_0)$ and a positive constant $C$, both depending on $\theta,\,R_0,\,M_\psi,\,M_A$ and $\eta_A$, such that for $\delta\in (0,\delta_0)$ and $\delta'\in (0,\theta\delta]$ there is a $C^k$ function $\ch$ defined on a neighborhood $\sv$ of $\overline\cy_\delta$ for which the following hold:
\medskip

\begin{itemize}
\item $0\le  \ch\le1$, $ \ch=0$ on $\sv-\scy_\delta$ and $ \ch=1$ on $\scy_{\delta'}$
\medskip

\item $|D^\alpha \ch(z)|\le C|\delta-\delta'|^{-|\alpha|}$ for $|\alpha|\le k$ and $z\in \sv$
\medskip

\item$\nabla \ch\cdot(A(z)N(z))=0$ for $z\in \partial\sw_\delta^+$.
\end{itemize}

\end{lemma}
\medskip

\begin{proof} Define a mapping $F:B_{R_0}\times\R\to\R^n$ by 
\begin{equation}\label{a.3} F(\tilde w,s)=(\tilde w,\psi(\tilde w))-s(AN)(\tilde w,\psi(\tilde w)).
\end{equation}
Then $F$ is $C^k$ and its derivative at the origin is invertible by \eqref{a.2}. Therefore for sufficiently small $\delta_0$  there is a neighborhood $S$ of the origin such that $F$ is a $C^k$ diffeomorphism of $S$ onto a neighborhood $\sv$ of $\ov\scy_{\delta_0}$. 

Now let $\delta,\delta'$ and $\theta$ be given as in the statement and define $\ep_j=\delta'+j(\delta-\delta')/3$ for $j=1,2$. Then choose a $C^\infty$ function $\tilde \ch_1$ on $\R^n$ which is a decreasing function of $|\tilde w|$ and whose value at $(\tilde w,s)$ is one or zero according as $|\tilde w|$ is less than $\ep_1$ or greater than $\ep_2$. We claim that if $\ch_1\equiv \tilde \ch_1\circ F^{-1}$, which is thus a $C^k$ function on $\sv$ taking values in $[0,1]$, then there is a constant $C_1$ depending on the same parameters as for $C$ in the statement of the theorem such that, if $z=(\tilde z,z_n)\in \sv$ with $|z_n|\le (\delta-\delta')/C_1$ and if $\delta_0$ is reduced further if necessary, then $ \ch_1(z)$ is one or zero according as $|\tilde z|$ is less than $\delta'$ or greater than $\delta$. The proof consists in checking that with $|z_n|$ so restricted and with $z=(\tilde z,z_n)=F(\tilde w,s)\in \sv$ then $|\tilde w|$ is less than $\ep_1$ or greater than $\ep_2$ according as $|\tilde z|$ is less than $\delta'$ or greater than $\delta$; details are left to the reader. 

Three more cutoff functions are required: Let $ \ch_2(\tilde z, z_n)$ be a decreasing $C^\infty$ function of $|\tilde z|$ which is one or zero according as $|\tilde z|$ is less than $\delta'$ or greater than $\delta$; let $\varphi$ be a decreasing $C^\infty$ function of $|z_n|$ which is one or zero according as $|z_n|$ is less than $(\delta-\delta')/12C_1$ or greater than $(\delta-\delta')/6C_1$; and let $\zeta$ be a decreasing $C^\infty$ function of $|z_n|$ which is one or zero according as $|z_n|$ is less than $\delta'$ or greater than $\delta$. We then take
$$ \ch(\tilde z,z_n) = \zeta(z_n)\big[\varphi(z_n) \ch_1(\tilde z, z_n)+(1-\varphi(z_n)) \ch_2(\tilde z,z_n)\big].$$
The first two items in the statement of the theorem are then easily checked. For the third we note that if $z=(\tilde z,z_n)\in \partial\sw_\delta^+$ then $|z_n|=|\psi(\tilde z)|\le C|\tilde z|^2$. Reducing $\delta_0$ further if necessary we find that $\varphi(z_n)=\zeta(z_n)=1$ and therefore $ \ch(z)= \ch_1(z)$, and in fact this holds as well for $z$ in a neighborhood of $\partial\sw_{\delta_0}^+$. Thus if $z=F(\tilde w,s)$ then  $ \ch(z)=\tilde \ch_1(\tilde w,s)$, which  is independent of $s$. Differentiating with respect to $s$ and recalling the definition \eqref{a.3}, we conclude that $\nabla \ch\cdot(AN) = 0$ at $z$ as required.
\medskip

\end{proof}

The following is the main result of this section:

\medskip

\begin{theorem} Assume that the hypotheses and notations in the first paragraph of this section are in force, let $\theta\in (0,1)$ be given and let $\delta_0\in (0,R_0)$ be as in {\rm Lemma~7.1}. Fix $\delta\in (0,\delta_0)$ and $\delta'\in (0,\theta\delta]$ and assume that $u\in H^1(\sw_\delta)$ satisfies one of the following:
\medskip 

\setlength{\hangindent}{18pt}
\noindent {\rm (D)} $u=0$ on $\partial\sw_\delta^+$ in the sense of traces and 
\begin{equation}\label{a.31}
\int_{\sw_\delta}(A\,\nabla u)\cdot\nabla v = F\cdot v
\end{equation}
for all $v\in H^1_0(\sw_\delta)$ and for some $F$ in the dual of $H^1_0(\sw_\delta)$.
\bigskip

\setlength{\hangindent}{18pt}
\noindent {\rm (N)} the relation \eqref{a.31} holds for all $v\in H^1(\sw_\delta)$ which vanish in the sense of traces on the sides and bottom of $\sw_\delta$ and for some $F$ in the dual of $H^1(\sw_\delta)$.
\medskip

\noindent Then the following hold:
\medskip

\setlength{\hangindent}{18pt}
\noindent {\rm (a)} There is a constant $C$ depending on $\theta,\,R_0,\,M_\psi(k),\,M_A(k)$ and $\eta_A$ such that if
$$|F\cdot v|\le a|v|_{L^2(\sw_\delta)} + b|\nabla v|_{L^2(\sw_\delta)} $$
for some constants $a$ and $b$ and for all $v$ as described in the respective cases {\rm (D)} or {\rm (N)} then
$$|\nabla u|^2_{L^2(\sw_{\delta'})} \le C\Big[(\delta-\delta')^{-2}|u|^2_{L^2(\sw_\delta)} + a|u|_{L^2(\sw_\delta)} +b^2\Big].$$

\medskip

\setlength{\hangindent}{18pt}
\noindent {\rm (b)} There is a constant $C$ depending on $k,\,\theta,\,R_0,\,M_\psi(k),\,M_A(k)$ and $\eta_A$ and on an upper bound for $(\delta-\delta')$ such that if 
$$F\cdot v=\displaystyle\int_{\sw_\delta} \big(\varphi\, v+\Phi\cdot\nabla v\big)$$
for some $\varphi\in H^{\min\{k-2,0\}}(\sw_\delta)$ and $\Phi\in [H^{k-1}(\sw_\delta)]^n$ then in either case {\rm (D)} or {\rm (N)} the solution $u$ is in $H^k(\sw_{\delta'})$ and for $j=1,\ldots,k$

\begin{align} \begin{split}\label{a.32}
\sum_{|\alpha|=j}|D^\alpha u|_{L^2(\sw_{\delta'})} \le C\ep^{-j}\Big[|u|_{L^2}+&\ep^{2}|\varphi|_{L^2}+\ep|\Phi|_{L^2}\\
	&+\sum_{i=1}^{j-1}\ep^{i+1}\big(|\varphi|_{H^{i-1}}+|\Phi|_{H^{i}}\big)\Big]
	\end{split}\end{align}
where  $\ep\equiv \delta-\delta'$, the norms on the right are taken over  $\sw_{\delta}$ and the sum on the right is omitted if $j=1$. 
\medskip 

\setlength{\hangindent}{18pt}
\noindent {\rm (c)} Apart from the hypotheses in {\rm (a)} and {\rm (b)} let $\theta'\in (0,\theta)$ be given and let $p\in (1,2]$ and $q= np/(n-2)$. Then there is a constant $C$ depending on $\theta,\,\theta',\, R_0,\,M_\psi(k)$, $M_A(k),\,\eta_A$ and $p$ such that if $\delta'\in[\theta'\delta,\theta\delta]$ and if $F\cdot v=\displaystyle\int_{\sw_\delta} \varphi\, v$ for some $\varphi\in L^p(\sw_\delta)$ then in either case {\rm (D)} or {\rm (N)} the solution $u$ is in $L^q(\sw_{\delta'})$ and 
\begin{equation}\label{a.312}
|u|_{L^q(\sw_{\delta'})}\le C\big[(\delta-\delta')^{-2/p}|u|_{L^p(\sw_{\delta})} + (\delta-\delta')^{2(p-1)/p}|\varphi|_{L^p(\sw_{\delta})}\big].
\end{equation}

\end{theorem}
\bigskip

\noindent{\bf Proof of Theorem 7.2(a).} Let $ \ch$ be as in Lemma 7.1. Then $v= \ch u$ is an allowable test function in \eqref{a.31} in either case (D) or (N) so that by \eqref{a.2} 
\begin{align*}
\int_{\sw_\delta'}|\nabla u&|^2\le C\int_{\sw_\delta} \ch (A\nabla u)\cdot\nabla u\\
	&=C\int_{\sw_\delta}\big[(A\nabla u)\cdot\nabla (\chi u)-u \nabla u\cdot (A^{tr}\nabla  \ch)\big]\\
	&=C\big[F\cdot ( \ch u) -\int_{\partial\sw_\delta^+} \textstyle \frac{1}{2}\displaystyle u^2\cdot(A^{tr}\nabla \ch)\cdot\nu+\int_{\sw_\delta}\textstyle \frac{1}{2}\displaystyle u^2\,{\rm div}(A^{tr}\nabla  \ch)\big]
\end{align*}
where $A^{tr}$ is the transpose of $A$. The boundary integral here is zero by the third item in Lemma 7.1, and bounds for the other two terms are obtained from the second item in Lemma 7.1 and the hypotheses on $F$.
\medskip

\noindent{\bf Proof of Theorem 7.2(b).}  This result is a special case of Proposition~4.3 of \cite{hoffbook} for much more general problems but whose proof is correspondingly more technically involved than is required here. This.proof is a standard argument in which $\sw_{\delta_0}$ is mapped to a region whose upper boundary is contained in a horizontal plane so that bounds for difference quotients of the transformed solution in the horizontal direction can be computed. These bounds are uniform in the increment, thus proving weak differentiability and providing bounds in $L^2$ for the corresponding horizontal derivatives; bounds for nonhorizontal derivatives are then obtained from \eqref{a.1} via the the ellipticity assumption \eqref{a.2}. We will omit the difference quotient argument, which is familiar to most readers, assuming therefore that $u\in H^k(\sw_{\delta_0})$, 
and content ourselves with deriving the bound in \eqref{a.32}, which is the point of interest here. 
 
Fix $\delta$ and $\delta'$ as in the statement and let $\delta'=\ep_k<\ep_{k-1}<\ldots <\ep_0=\delta$ so that $\ep_j-\ep_{j+1}=\ep/k$. The hypothesis of (a) then holds with $\delta,\delta'$ replaced by $\ep_0,\ep_1$ and with $a=|\varphi|_{L^2(\sw_{\delta})}$ and $b=|\Phi|_{L^2(\sw_{\delta})}$. Thus
\begin{equation}\label{a.50}
\int_{\sw_{\ep_1}}|\nabla u|^2 \le C\big(\ep^{-2}|u|_{L^2(\sw_{\delta})}^2+\ep^2   |\varphi|^2_{L^2(\sw_{\delta})}+|\Phi|^2_{L^2(\sw_{\delta})}\big)\equiv M_1. 
\end{equation}
This proves \eqref{a.32} for $k=1$. We therefore assume that $k\ge 2$ and
flatten $\partial\sw_{\delta_0}^+$ as follows. Let $Y$ be the $C^{k+1}$ map defined by $Y(\tilde z,z_n)=(\tilde y,y_n)$ where $\tilde y=\tilde z$ and $y_n=z_n-\psi(\tilde z)$ and let $\sw'_{\delta_0}=Y(\sw_{\delta_0})$.  Define $u'\in H^k(\sw'_{\delta_0})$ and $\varphi'\in H^{k-2}(\sw'_{\delta_0})$ by $u=u'\circ Y$ and $\varphi=\varphi'\circ Y$, and for allowable test functions $v$ let $v=v'\circ Y$. It then follows that
\begin{equation}\label{a.51}
\int_{\sw'_{\delta_0}}(A'\,\nabla u')\cdot\nabla v'= \int_{\sw'_{\delta_0}} \big(\varphi'\,v'+\Phi'\cdot\nabla v'\big)
\end{equation}
where 
\begin{equation}\label{a.52}
A'=SAS^{tr}, \ S=\partial Y/\partial z\ {\rm and}\ \Phi'=S(\Phi\circ Y ).
\end{equation} 
Also, if we define
$$\tilde M_j=\sum_{1\le|\alpha|\le j}\int_{\sw'_{\ep_j}}|D^\alpha u'|^2$$
where $\sw'_{\ep_j}=Y(\sw_{\ep_j})$, then by \eqref{a.50} 
\begin{equation}\label{a.53} \tilde M_1\le C\big(\ep^{-2}|u|_{L^2(\sw_{\delta})}^2+\ep^2|\varphi|^2_{L^2(\sw_{\delta})}+|\Phi|^2_{L^2(\sw_{\delta})}\big).
\end{equation}
\medskip
  
Now suppose that $v\in H^{k}(\sw_{\delta_0})$ is an allowable test function which is zero on neighborhoods of the sides and bottom of $\sw_{\delta_0}$ and let $v'$ be as above. Then if $0<|\alpha|\le k-1$ and $D_y^\alpha$ is a horizontal derivative, that is, $\alpha=(\alpha_1,\ldots,\alpha_n)$ with $\alpha_n=0$, then $D_y^\alpha v'$ is an allowable test function in \eqref{a.51}:
\begin{equation*}
\int_{\sw'_{\delta_0}}(A'\,\nabla u')\cdot\nabla (D^\alpha v')= \int_{\sw'_{\delta_0}} \big( \varphi'\,D^\alpha v'+\Phi'\cdot\nabla D^\alpha v'\big).
\end{equation*}
We integrate by parts on the left, noting that the boundary integrals are zero, and apply the Leibnitz formula to bound $D^\alpha (A'\nabla u')$. On the right side we let $\alpha=\beta+\gamma$ where $|\gamma|=1$ and $|\beta|\le k-2$, then integrate by parts again. The result is that
\begin{align}\begin{split}\label{a.54}
\int_{\sw'_{\delta_0}}(A'\,\nabla (D^\alpha u'))&\cdot\nabla v'\\
= (-1)^{|\alpha|+|\beta|}&\int_{\sw'_{\delta_0}} \big[(D^\beta \varphi')\,(D^\gamma v')+(D^\beta\Phi')\cdot\nabla (D^\gamma v')\big]\\
&\qquad\qquad	+O\Big(\sum_{|\eta|\le |\alpha|}\int_{\sw'_{\delta_0}} |D^\eta u'|\,|\nabla v'| \Big) 
	\end{split}
\end{align}
(the second term on the right is a generic quantity bounded by $C$ times the enclosed sum). A standard density argument shows that this holds as well for allowable test functions $v'=v\circ Y^{-1}\in H^1(\sw'_{\delta_0})$.

Next we fix $j\in \{1,\ldots,k-1\}$ and derive a bound for $\tilde M_{j+1}$ in terms of $\tilde M_j$. Let $\chi$ be as in Theorem~7.1 with $\delta, \delta'$ replaced by $\ep_j,\ep_{j+1}$. Define $\chi'$ by $\chi=\chi'\circ Y$ and let $\alpha$ be a horizontal derivative as above with $|\alpha|=j$. Then since $u'\in H^k(\sw'_{\delta_0})$, \eqref{a.54} holds with $v'=\chi'D^\alpha u'$:
\begin{align}\begin{split}\label{a.55}
\int_{\sw'_{\ep_j}}(A'\,\nabla (D^\alpha u'))&\cdot\nabla (\chi'D^\alpha u')\\
= (-1)^{|\alpha|+|\beta|}&\int_{\sw'_{\ep_j}} \big[(D^\beta \varphi')\,(D^\gamma (\chi'D^\alpha u'))+(D^\beta\Phi')\cdot \nabla D^\gamma(\chi'D^\alpha u')\big]\\
	&\qquad\qquad\qquad+O\big(\sum_{|\eta|\le j}\int_{\sw'_{\ep_j}} |D^\eta u'|\,|\nabla (\chi' D^\alpha u')| \big). 
	\end{split}
\end{align}
The second item in Theorem 7.1 shows that the first two terms on the right here are bounded by 
 $$C\big(|\varphi'|_{{H^{j-1}(\sw'_{\ep_j})}}+|\Phi'|_{{H^{j}(\sw'_{\ep_j})}}\big)\Big[\Big(\int_{\sw'_{\ep_j}}\chi'|\nabla D^\alpha u'|^2\Big)^{1/2} + \ep^{-1}\tilde M_j^{1/2}\Big]$$ 
and the third is bounded by
$$C\tilde M_j^{1/2}\Big[\Big(\int_{\sw'_{\ep_j}}\chi'|\nabla D^\alpha u'|^2\Big)^{1/2}+\ep^{-1}\Big].$$
The term on the left in \eqref{a.55} is
\begin{align*}
\int_{\sw'_{\ep_j}}(A'\,\nabla &(D^\alpha u'))\cdot\big[\chi' \nabla(D^\alpha u')+(D^\alpha u')\nabla\chi'\big]\\
	&\ge \eta_A\int_{\sw'_{\ep_j}}\chi'|\nabla (D^\alpha u')|^2+{\textstyle\frac{1}{2}}\int_{\sw'_{\ep_j}}(A'\nabla(D^\alpha{u'}^2))\cdot\nabla\chi'
\end{align*}
by \eqref{a.2}. The first term on the right here absorbs the integrals in the previous two displays. For the second term we again integrate by parts, noting that the integrand is zero on the bottom and sides of the domain and that the normal on the top is the standard basis vector $e_n\equiv(0,\ldots,1)$. The boundary integral in this integration by parts is therefore zero because by direct computation from the definition of $Y$, $(A'e_n)\cdot\nabla \chi'=(A\nu)\cdot\nabla\chi$, which is zero by the third item in Lemma~7.1. The resulting integral over $\sw'_{\ep_j}$  is then easily seen to be bounded by $C\ep^{-2}\tilde M_j$.  Assembling these three bounds and simplifying we obtain finally

\begin{equation*}
\int_{\sw'_{\ep_j}}\chi'|\nabla (D^\alpha u')|^2\le C\big(\ep^{-2}\tilde M_j+|\varphi'|^2_{H^{j-1}(\sw'_{\ep_j})}+|\Phi'|^2_{H^{j}(\sw'_{\ep_j})}).
\end{equation*}
Recall that this holds for horizontal derivatives $D^\alpha$ of order $\nobreak{j\le k-1}$. We now change notation and let $\alpha$ denote a multi-index of length $\nobreak{j+1\le k}$. Then since $\chi'\equiv 1$ on $\sw'_{\ep_{j+1}}$, 
\begin{equation}\label{a.56}
\int_{\sw'_{\ep_{j+1}}}|D^\alpha u'|^2\le C\big(\ep^{-2}\tilde M_j+|\varphi'|^2_{H^{j-1}(\sw'_{\ep_j}) }+|\Phi'|^2_{H^{j}(\sw'_{\ep_j})}\big)
\end{equation}
for $|\alpha|\le j+1$ provided that $\alpha_n$ is zero or one. 

Now consider the case $j=1$. By taking test functions $v'$ of compact support in $\sw'_{\delta_0}$ in \eqref{a.51} we find that 
\begin{equation}\label{a.57}
{\rm div}\,(A'\nabla u')=\varphi'
\end{equation}
in the weak sense on $\sw'_{\delta_0}$. Also, if $a^{i,j}$ are the components of the matrix $A'$ then \eqref{a.56} shows that the terms $a^{i,j}u'_{y_i,y_j}$ in \eqref{a.57} for $(i,j)\not=(n,n)$ are in $L^2(\sw'_{\ep_{2}})$ with the squares of their norms bounded by the right side of \eqref{a.56} with $j=2$, as are $\varphi'$ and the lower order terms;  hence so is the term $a^{n,n}u'_{y_n,y_n}$ (with a new constant $C$). Putting $\xi=e_n$ in \eqref{a.2} and recalling that $A'=SAS^{tr}$, we find that $a^{n,n}\ge C^{-1}$. Thus \eqref{a.56} holds for {\it all} $\alpha$ of length two, that is,   
$$\tilde M_2\le C\big(\ep^{-2}\tilde M_1+|\varphi'|^2_{L^2(\sw'_{\ep_1})}+|\Phi'|^2_{H^{1}(\sw'_{\ep_1})}\big).$$
More generally, for $j\le k-1$ and $\alpha$ of length $j+1$ we write $\alpha=\beta + \gamma$ with $\gamma=2$, differentiate \eqref{a.57} $j-1$ times with successively larger values of $\alpha_n$ and apply \eqref{a.56} to find that \eqref{a.56} holds for all $\alpha$ of length $j+1$. The conclusion is that
\begin{equation}\label{a.58}
\tilde M_{j+1}\le C\big(\ep^{-2}\tilde M_j+|\varphi'|^2_{H^{j-1}(\sw'_{\ep_j})}+|\Phi'|^2_{H^{j}(\sw'_{\ep_j})}\big)
\end{equation}
for $j=1,\ldots,k-1$ and for a constant $C$ as described in the statement of the theorem. To complete the proof we proceed by induction on $j$ starting from \eqref{a.53}, then convert the resulting bounds for $\tilde M_j$ by way of the $C^{k+1}$ diffeomorphism $Y$ to obtain the bound in \eqref{a.32}. 
\medskip

\noindent{\bf Proof of Theorem 7.2(c).}  The idea of the proof is to apply \eqref{a.31} with test function $v= \ch |u|^{p/2}$ where $ \ch$ is as in Lemma 7.1. Issues of regularity and integrability arise near points where $u=0$, however. To avoid these we introduce $C^\infty$ functions on $\R$ defined for $\ep>0$ by 
\begin{itemize}
\item $f_\ep(t)=\displaystyle\frac{t^2}{\sqrt{t^2+\ep^2}}+\ep$
\medskip
\item $g_\ep'=f_\ep^{p-2}(f_\ep')^2$ and $h'_\ep=g_\ep$
\medskip
\item $g_\ep(0)=h_\ep(0)=0$.
\end{itemize}
In particular, $f_\ep$ is a $C^\infty$ approximation to the absolute value function. Elementary arguments prove the following, in which we denote by $C'$ a universal positive constant independent of present considerations and by $C'_\ep$ a similar constant but which may also depend on $\ep$:
\medskip

\begin{itemize}
\item $\ep\le f_\ep(t)\le |t|+\ep\ ,|f'_\ep(t)|\le C'$ and $f_\ep(t)\to |t|$ uniformly on $\R$ as $\ep\to 0$
\medskip

\item $|g_\ep(t)|\le C'(|t|^{p-1}+\ep^{p-1})$ and $|g_\ep'(t)|\le C'_\ep$
\medskip

\item $|h_\ep(t)|\le C'(|t|^p+\ep^{p-1}|t|)$.
\end{itemize}
\medskip
We denote $f_\ep\circ u,\, g_\ep\circ u$ and $h_\ep\circ u$ by $F_\ep,\,G_\ep$ and $H_\ep$ and we let $ \ch$ be as in Lemma~7.1. It then follows from the above that $F^{p/2}_\ep\in H^1(\sw_\delta)$ and that $G_\ep$  is an allowable test function in \eqref{a.31}. Thus
\begin{align*}
\int_{\sw_{\delta'}}|\nabla (F_\ep^{p/2})|^2\le &C\int_{\sw_\delta}  \ch(A\nabla(F_\ep^{p/2}))\cdot\nabla(F_\ep^{p/2})=C\int_{\sw_\delta}  \ch(A\nabla u)\cdot\nabla G_\ep\\
	&=C\Big[\int_{\sw_\delta}(A\nabla u)\cdot\nabla ( \ch G_\ep)-\nabla H_\ep\cdot (A^{tr}\nabla \ch)\Big] \\
		&=C\Big[\int_{\sw_\delta}\varphi  \ch G_\ep -\int_{\partial\sw^+_\delta} H_\ep(A^{tr}\nabla \ch)\cdot\nu \\&\qquad\qquad\qquad\qquad\qquad+	\int_{\sw_\delta} H_\ep\,{\rm div}(A^{tr}\nabla \ch)\Big].	
\end{align*}
The boundary integral here is zero by the third item in Lemma 7.1 and bounds for the other two terms on the right are obtained from the second item in Lemma 7.1, the assumed bounds for $A$ and the above properties of $g_\ep$ and $h_\ep$. These together with H\"older's inequality thus show that
\begin{align}\begin{split}\label{a.41}|\int_{\sw_{\delta'}}|\nabla (F_\ep^{p/2})|^2\le C&\Big[(\delta-\delta')^{2(p-1)}|\varphi|^p_{L^p}+(\delta-\delta')^{-2}|u|_{L^p}^p\Big]\\ +&C\ep^{p-1}\Big[|\varphi|_{L^1}+[(\delta-\delta')^{-2}|u|_{L^1}\Big]
\end{split}\end{align}
where the norms on the right are taken over $\sw_\delta$. 
The scaling result in Theorem~7.3(b) below (whose proof is completely independent of all other results in this paper) shows that there is a constant $C$ independent of $\delta'$ such that for $v\in H^1(\sw_{\delta'})$, 
$$|v|_{L^{2n/(n-2)}}\le C\big[(\delta')^{-1}|v|_{L^2}+|\nabla v|_{L^2}\big]$$
(norms are over $\sw_{\delta'}$). Applying this together with the obvious bound $|F_\ep^{p/2}|_{L^2}\le C(|u|^{p/2}_{L^p}+\ep^{p/2})$ and noting that $\delta'\ge 
\theta_1(\delta-\delta')/(1-\theta_2)$ we find that 
$$
\limsup_{\ep\to 0}|F_\ep^{p/2}|_{L^{2n/(n-2)}(\sw_{\delta'})} \le C\big[(\delta-\delta')^{-1}|u|^{p/2}_{L^p(\sw_\delta)}+(\delta-\delta')^{p-1}|\varphi|^{p/2}_{L^p(\sw_\delta)}\big].
$$
The monotone convergence theorem shows that the limit on the left is $|u|^{2/p}_{L^q(\sw_{\delta'})}$, and this proves the bound in \eqref{a.312}. 

\rightline\qedsymbol

\medskip

In the following theorem we consider the problem \eqref{a.1} posed in a ball $B_{R_0}$ and for test functions in $H^1_0(B_{R_0})$. Essentially the same results as in Theorem 7.2 hold, but the proofs are greatly simplified because boundary considerations do not enter. These proofs are therefore omitted.
	
 \medskip

\begin{theorem} Let $n\ge 2$ and fix $R_0>0$ and an $n\times n$-matrix valued function $A$ on $B_{R_0}\subset\R^n$ whose derivatives up to order $k\ge 1$ exist and are bounded in $B_{R_0}$ by a constant $M_A$ and which satisfies the positive definiteness condition \eqref{a.2} on $B_{R_0}$ with constant $\eta_A$. Assume that $u\in H^1(B_{R_0})$ satisfies \eqref{a.1} with $\sw_{R_0}$ replaced by $B_{R_0}$, for all $v\in H^1_0(B_{R_0})$ and for some $F\in (H^1_0(B_{R_0}))^*$. Let $\theta\in (0,1)$ be given and fix $\delta\in (0,R_0)$ and $\delta'\in (0,\theta\delta]$. Then all the statements in parts {\rm (a), (b) and (c)} of {\rm Theorem 7.2} hold with $\sw_\delta$ and $\sw_{\delta'}$ replaced by $B_\delta$ and $B_{\delta'}$ and with references to $M_\psi$ omitted.
\end{theorem}

In the final result of this section we derive scaled versions of the standard imbedding $H^1\mapsto L^{2n/(n-2)}$, which holds on open sets satisfying a ``cone condition," and $H^{k_0}\mapsto C^{0,\lambda_0}$, which holds on sets having a Lipschitz boundary (see \cite{adams} pg. 85 for precise definitions and statements). 

 \begin{theorem}

 \setlength{\hangindent}{18pt}
\noindent {\rm (a)}  Let $\sw_\delta\subset\R^n$, $n\ge 2$, be as in the first paragraph of this section with $\delta\le R_0$ but with the exception that $\psi$ is assumed only to be Lipschitz continuous with constant $L$; and let  $k_0$ and $\lambda_0$ be as in \eqref{2.9}. Then there are constants $C(n,R_0,L)$ and $C(n,R_0,L,\lambda_0)$ such that if $v\in H^{k_0}(\sw_\delta)$ then $v$ is H\"older continuous with exponent $\lambda_0$ (modulo equivalence class) and  satisfies
\begin{equation}\label{a.61}
\sup_{x\in\sw_\delta}|v(x)|\le C(n,R_0,L)\sum_{j=0}^{k_0}\sum_{|\alpha|=j}\delta^{j-n/2}|D_x^\alpha v|_{L^2(\sw_\delta)}
\end{equation} and
\begin{equation}\label{a.62}
\langle v\rangle_{\sw_\delta,\lambda_0}\le C(n,R_0,L,\lambda_0)\sum_{j=0}^{k_0}\sum_{|\alpha|=j}\delta^{j-n/2-\lambda_0}|D_x^\alpha v|_{L^2(\sw_\delta)}.
\end{equation}
\medskip

\setlength{\hangindent}{18pt}
\noindent {\rm (b)} Let $\sv$ be an open set in $\R^n$, $n\ge 2$, satisfying a $K_\delta$ cone condition; that is, there is a subset $S$ of the $n-1$-sphere in $\R^n$ with positive surface measure $|S|$ and a positive number $\delta>0$ such that if 
$$K_\delta=\{r\omega : 0\le r\le \delta\ {\rm and\ }\omega\in S\}$$ 
then every point in $\sv$ is the vertex of a cone contained in $\sv$ and congruent to $K_\delta$. There is then a constant $C(n,|S|)$ independent of $\delta$ such that if $v\in H^1(\sv)$ then $v\in L^{2n/(n-2)}(\sv)$ and
$$|v|_{L^{2n/(n-2)}(\sv)}\le C\big(\delta^{-1}|v|_{L^{2}(\sv)} + |\nabla v|_{L^{2}(\sv)}\big).$$

\end{theorem}
\medskip

{\noindent\bf Proof.} To prove (a) we let $Y(x)=R_0x/\delta$ and define
	$$\sw'\equiv Y(\sw_\delta)=\{(\tilde y,y_n): |\tilde y|< R_0\ {\rm and}\ -R_0<y_n<\varphi(\tilde y)\}$$
where $\varphi= R_0(\psi\circ Y^{-1})/\delta$. Then $\varphi$ is Lipschitz with constant $L$ and the embedding $H^{k_0}(\sw')\mapsto C^{0,\lambda_0}(W')$ therefore holds with constant $C(n,R_0,L)$. Given $v\in H^{k_0}(\sw_\delta)$ we then apply this embedding to the function $w\equiv v\circ Y^{-1}$, noting that
\begin{align*}
|D_y^\alpha w|^2_{L^2(\sw')}=\int_{\sw_\delta}(\delta/R_0)^{2|\alpha|} |D_x^\alpha v(x)|^2(R_0/\delta)^n dx\\=C(R_0)\delta^{2|\alpha|-n}|v|^2_{L^2(\sw_\delta)},
\end{align*}
then rescale to obtain \eqref{a.61} and \eqref{a.62}.

To prove (b) we define $\sv_1\equiv\{x/\delta : x\in \sv\}$. Then $\sv_1$ satisfies a $K_1$-cone condition so that the embedding $H^1(\sv_1)\to L^{2n/(n-2)}(\sv_1)$ holds with a constant $C$ which is independent of $\delta$. Given $v\in H^1(\sv)$ we then apply this embedding to the element $w$ defined by $w(y)=v(\delta y)$, then rescale to obtain the required bound for $v$.
\medskip

\rightline\qedsymbol

\bibliographystyle{amsplain}
\bibliography{citations}
\smallskip

\end{document}